\documentclass [11pt]{amsart}

\usepackage{graphicx}
\usepackage{amsmath}
\usepackage{amssymb}
\usepackage{amsthm}

\newcommand{\nin}{\notin}
\newcommand{\bdry}{\partial}

\renewcommand{\Tilde}{\widetilde}

\newcommand{\cardioid}{\heartsuit}
\newcommand{\Fibo}{\mathbb{F}\hspace{-.2pt}\text{\rm \i bo}}

\newtheorem{thm}{Theorem}[section]
\newtheorem{prop}[thm]{Proposition}
\newtheorem{lemma}[thm]{Lemma}
\newtheorem{corol}[thm]{Corollary}
\newtheorem*{conj}{Conjecture}
\newtheorem*{thm*}{Theorem}

\newtheorem*{Lthm}{Lyubich's Theorem}

\theoremstyle{definition}
\newtheorem{defn}{Definition}[section]

\theoremstyle{remark}

\newtheorem*{note}{Note}

\begin{document}

\title{Geometry of $Q$-recurrent maps}
\author{Rodrigo A. P\'erez}
\address{Department of Mathematics, Cornell University, Ithaca, NY
14853. USA.}
\email{rperez@math.cornell.edu}
\thanks{Research supported by an NSF Postdoctoral Fellowship in the
Mathematical Sciences, grant DMS-0202519.}

\begin{abstract}
  {\it Given a critically periodic quadratic map with no secondary
  renormalizations, we introduce the notion of $Q$-recurrent quadratic
  polynomials. We show that the pieces of the principal nest of a
  $Q$-recurrent map $f_c$ converge in shape to the Julia set of $Q$. We use
  this fact to compute analytic invariants of the nest of $f_c$, to give a
  complete characterization of complex quadratic Fibonacci maps and to obtain
  a new auto-similarity result on the Mandelbrot set.}
\end{abstract}

\maketitle

\section{Introduction}\label{sect:Intro}
The principal nest of a quadratic polynomial is a collection of pieces of the
puzzle with good shrinking properties (see \cite{L_nest}). In \cite{1st_part}
we proposed a supplementary construction, the {\it frame} of the nest, to help
classify different nest types. Based on that classification, we introduce in
this paper the family of $Q$-recurrent maps. In essence, these are polynomials
for which the nest frames are combinatorially modeled on the puzzle of a
chosen critically periodic quadratic polynomial $Q$. We will generalize to the
family of $Q$-recurrent maps some results of M. Lyubich and L. Wenstrom that
concern the real Fibonacci parameter $c_{\rm fib}$. It is shown in
\cite{L_teichmuller} that the central pieces of the nest of $f_{c_{\rm fib}}$
are asymptotically similar to the Julia set of $z^2-1$. This fact is expanded
in \cite{W}, where the shape of pieces is exploited to compute the exact rate
of growth of the principal nest moduli. Wenstrom translates these results to
the parameter plane in order to obtain the corresponding shape of parapieces
and the rate of paramoduli growth around $c_{\rm fib}$ on the Mandelbrot set
$M$.

The corresponding results for $Q$-recurrent maps are as follows:

\begin{thm*}
  Let $c$ be the center of a prime component of $M$ and $Q$ its associated
  quadratic polynomial. Then the nest pieces of a $Q$-recurrent map converge
  in shape to the Julia set $K_Q$.
\end{thm*}

The asymptotic shape of pieces is used to derive analytic information about
the nest. In the simplest case, $Q(z)=z^2-1$, the family of
$(z^2-1)$-recurrent maps has an interesting description.

\begin{thm*}
  A quadratic polynomial is $(z^2-1)$-recurrent if and only if it is complex
  Fibonacci. Moreover, the principal moduli of any map in this family grow
  linearly with growth rate $\frac{\ln 2}3$.
\end{thm*}

The corresponding result is different when the critical point of $Q$ has
period $\geq 3$.

\begin{thm*}
  For $Q(z) \neq z^2-1$ the principal moduli of a $Q$-recurrent map grow
  exponentially with a rate that depends on the period of $Q$.
\end{thm*}

When the shape results are translated to the Mandelbrot set, the above
statements find parametric counterparts.

\begin{thm*}
  The paranest pieces around a $Q$-recurrent parameter $c$ converge in shape
  to $K_Q$. The growth of paramoduli around $c$ is as in the corresponding
  dynamical plane. The family of $Q$-recurrent parameters forms a dyadic
  Cantor set of Hausdorff dimension 0.
\end{thm*}

As an application of parapiece shape, an easy diagonal argument yields a
powerful auto-similarity result on $M$.

\begin{thm*}
  Let $c_1, c_2 \in \bdry M$ be two parameters such that $f_{c_2}$ has no
  indifferent periodic orbits that are rational or linearizable. Then there
  exists a sequence of parapieces $\{ \Upsilon_1, \Upsilon_2, \ldots \}$
  converging to $c_1$ as compact sets, but such that $\Upsilon_n$ converges to
  $K_{c_2}$ in shape.
\end{thm*}

\subsection{Paper structure}\label{subsect:Contents}

The puzzle of Yoccoz, the principal nest and their parametric counterparts are
defined in Section \ref{sect:Basics} as a means to introduce notation. The
adjacency graphs introduced in Subsection \ref{subsect:Graphs} are used in
Section \ref{sect:Frames} to define the frame of a principal nest, and in
Section \ref{sect:Def_Q-rec} to describe $Q$-recurrent behavior.

In Section \ref{sect:Q_Recurrency} we present the classification of complex
Fibonacci quadratic polynomials and the results on the shape of pieces and
growth of moduli for the nest of $Q$-recurrent maps.

The corresponding results on parametric pieces of the Mandelbrot set $M$ are
presented in Section \ref{sect:Parameter_Space}. Theorem
\ref{thm:M_Similarity} introduces a new similarity phenomenon between
different locations of $\bdry M$.

An appendix summarizes the tools borrowed from complex variables. This
includes an extension of the Gr\"otzsch inequality and brief discussions of
Carath\'eodory topology, Koebe's distortion lemma and the Teichm\"uller space
of a surface.

\subsection{Acknowledgments}
Many thanks are due to Mikhail Lyubich and John Milnor for their helpful
suggestions. Some of the pictures were created with the PC program {\tt
mandel.exe} by Wolf Jung \cite{J}.

\section{Basic notions}\label{sect:Basics}
\subsection{Basic complex dynamics}
In order to fix notation, let us start by defining the basic notions of
complex dynamics that will be used; the reader is referred to \cite{DH_orsay}
and \cite{M_book} for details on this introductory material.

We focus attention on the {\it quadratic family} $\mathcal{Q} := \big\{ f_c:z
\mapsto z^2 + c \mid c \in \mathbb{C} \big\}$. For every $c$, the compact sets
$K_c := \big\{ z \mid \text{the sequence } \{f_c^{\circ n}(z)\} \text{ is
bounded} \big\}$ and $J_c := \bdry K_c$ are called the {\bf filled Julia set}
and {\bf Julia set} respectively. Depending on whether the orbit of the
critical point 0 is bounded or not, $J_c$ and $K_c$ are connected or totally
disconnected. The {\bf Mandelbrot set} is defined as $M := \big\{c \mid c
\in K_c \big\}$; that is, the set of parameters with bounded critical orbit.

A component of ${\rm int} \,M$ that contains a parameter with an attracting
periodic orbit will be called a {\bf hyperbolic component}\footnote{Though, of
course, it is conjectured that all interior components are hyperbolic.}. The
boundary of a hyperbolic component can either be real analytic, or fail to be
so at one cusp point. The later kind are called {\bf primitive} components. In
particular, the hyperbolic component $\cardioid$ associated to $z \mapsto z^2$
is bounded by a cardioid known as the {\bf main cardioid}.

$M$ contains infinitely many small homeomorphic copies of itself, accumulating
densely around $\bdry M$. In fact, every hyperbolic component $H$ other than
the main one is the base of one such small copy $M'$. $H$ is called {\bf
prime} if it is not contained in any other small copy. To simplify later
statements, prime components are further subdivided in {\bf immediate}
(non-primitive components that share a boundary point with $\cardioid$) and
{\bf maximal} (primitive components away from $\bdry\cardioid$).

\subsection{External rays, wakes and limbs}
  \label{subsect:External_Rays_n_Wakes} Since $f_c^{-1}(\infty) =\{\infty\}$,
the point $\infty$ is a fixed critical point and a classical result of
B\"ottcher yields a change of coordinates that conjugates $f_c$ to $z \mapsto
z^2$ in a neighborhood of $\infty$. With the requirement that the derivative
at $\infty$ is $1$, this conjugating map is denoted $\varphi_c:N_c
\longrightarrow \overline{\mathbb{C}} \setminus \overline{\mathbb{D}_R}$,
where $\mathbb{D}_R$ is the disk of radius $R \geq 1$ and $N_c$ is the maximal
domain of unimodality for $\varphi_c$. It can be shown that $N_c =
\overline{\mathbb{C}} \setminus K_c$ and $R=1$ whenever $c \in M$. Otherwise,
$N_c$ is the exterior of a figure $8$ curve that is real analytic and
symmetric with respect to 0. In this case, $R>1$ and $K_c$ is contained in the
two bounded regions determined by the $8$ curve.

Consider the system of radial lines and concentric circles in $\mathbb{C}
\setminus \mathbb{D}_R$ that characterizes polar coordinates. The pull-back of
these curves by $\varphi_c$, creates a collection of {\bf external rays}
$r_{\theta}$ $\big( \theta \in [0, 1) \big)$ and {\bf equipotential curves}
$e_s$ \big(here $s \in (R,\infty)$ is called the {\bf radius} of $e_s$\big) on
$N_c$. These form two orthogonal foliations that behave nicely under dynamics:
$f_c(r_{\theta}) = r_{2\theta}$, $f_c(e_s) = e_{(s^2)}$. When $c \in M$, we
say that a ray $r_{\theta}$ {\bf lands} at $z \in J_c$ if $z$ is the only
point of accumulation of $r_{\theta}$ on $J_c$. \\

A similar coordinate system exists around the Mandelbrot set. For $c \nin M$,
define the map
\begin{equation}\label{eqn:Phi_M}
  \Phi_M(c) := \varphi_c(c).
\end{equation}

In \cite{DH_orsay} it is shown that $\Phi_M:\mathbb{\overline{C}} \setminus M
\longrightarrow \mathbb{\overline{C}} \setminus \overline{\mathbb{D}}$ is a
conformal homeomorphism tangent to the identity at $\infty$. This yields
connectivity of $M$ and allows us to define {\bf parametric external rays} and
{\bf parametric equipotentials} as in the dynamical case. Since there is
little risk of confusion, we will use the same notation ($r_{\theta}, e_s$) to
denote these curves and say that a parametric ray lands at a point $c \in
\bdry M$ if $c$ is the only point of accumulation of the ray on $M$.

For the rest of this work, all rays considered, whether in dynamical or
parameter plane, will have rational angles. These are enough to work out our
combinatorial constructions and satisfy rather neat properties.

\begin{prop} (\cite{M_book}, ch.18) Both in the parametric and the dynamical
  situations, if $\theta \in \mathbb{Q}$ the external ray $r_{\theta}$
  lands. In the dynamical case, the landing point is (pre-)periodic with the
  period and preperiod determined by the binary expansion of $\theta$. A point
  in $J_c$ (respectively $\bdry M$) can be the landing point of at most, a
  finite number of rays (respectively parametric rays). If this number is
  larger than $1$, each component of the plane split by the landing rays will
  intersect $J_c$ (respectively $\bdry M$).
\end{prop}

Unless $c = \frac14$, $f_c$ has two distinct fixed points. If $c \in M$,
these can be distinguished since one of them is always the landing point of
the ray $r_0$. We call this fixed point $\beta$. The second fixed point is
called $\alpha$ and can be attracting, indifferent or repelling, depending on
whether the parameter $c$ belongs to $\cardioid$, $\bdry\cardioid$, or
$\mathbb{C} \setminus \overline{\cardioid}$. The map $\psi_0: \cardioid
\longrightarrow \mathbb{D}$ given by $c \mapsto f'_c(\alpha_c)$ is the Riemann
map of $\cardioid$ normalized by $\psi_0(0)=0$ and $\psi'_0(0)>0$. Since the
cardioid is a real analytic curve except at $\frac14$, $\psi_0$ extends to
$\overline{\cardioid}$.

The fixed point $\alpha$ is parabolic exactly at parameters $c_{\eta} \in
\bdry \cardioid$ of the form $c_{\eta} = \psi_0^{-1} \left( e^{2\pi i \eta}
\right)$ where $\eta \in \mathbb{Q} \cap [0,1)$. If $\eta \neq 0$, $c_{\eta}$
is the landing point of two parametric rays at angles $t^-(\eta)$ and
$t^+(\eta)$.

\begin{defn}
  The closure of the component of $\mathbb{C} \setminus \left( r_{t^-(\eta)}
  \cup c_{\eta} \cup r_{t^+(\eta)} \right)$ that does not contain $\cardioid$
  is called the {\bf $\eta$-wake} of $M$ and is denoted $W_{\eta}$. The {\bf
  $\eta$-limb} is defined as $L_{\eta} = M \cap W_{\eta}$.
\end{defn}

\begin{defn}
  Say that $\eta=\frac pq$, written in lowest terms. Then $\mathcal{P}\big(
  \frac pq \big)$ will denote the unique set of angles whose behavior under
  doubling is a cyclic permutation with combinatorial rotation number $\frac
  pq$.
\end{defn}

If $\mathcal{P}\big( \frac pq \big) = \{t_1,\ldots , t_q \}$, then for any
parameter $c \in L_{p/q}$ the corresponding point $\alpha$ splits $K_c$ in $q$
parts, separated by the $q$ rays $\{ r_{t_1}, \ldots , r_{t_q} \}$ landing at
$\alpha$. The two rays whose angles span the shortest arc separate the
critical point 0 from the critical value $c$. These two angles turn out to be
$t^-(\frac pq)$ and $t^+(\frac pq)$.

\subsection{Yoccoz puzzles}
The Yoccoz {\bf puzzle} is well defined for parameters $c \in L_{p/q}$ for any
any $\frac pq \in \mathbb{Q} \cap [0,1)$ with $(p,q)=1$. If 0 is not a
preimage of $\alpha$, the puzzle is defined at infinitely many depths and we
will restrict attention to these parameters. Since we describe properties of a
general parameter, it is best to omit subscripts and write $f$ instead of
$f_c$, $K$ instead of $K_c$ and so on.

Let us fix the neighborhood $U$ of $K$ bounded by the equipotential of radius
2. The rays that land at $\alpha$ determine a partition of $U \setminus \{
r_{t_1}, \ldots , r_{t_q} \}$ in $q$ connected components.  We will call the
closures $Y_0^{(0)}, Y_1^{(0)}, \ldots, Y_{q-1}^{(0)}$ of these components,
{\bf puzzle pieces} of depth $0$. At this stage the labeling is chosen so that
$0 \in Y_0^{(0)}$ and $f \left( K \cap Y_j^{(0)} \right) = K \cap
Y_{j+1}^{(0)}$; where the subindices are understood as residues modulo $q$. In
particular, $Y_1^{(0)}$ contains the critical value $c$ and the angles of its
bounding rays turn out to be $t^-(\frac pq),t^+(\frac pq)$.

The puzzle pieces $Y_i^{(n)}$ of higher depths are recursively defined as the
closures of every connected component in $f^{\circ (-n)} \left( \bigcup
{\rm int}\, Y_j^{(0)} \right)$; see Figure \ref{fig:Puzzle_Graphs}. At each
depth $n$, there is a unique piece which contains the critical point and we
will always choose the indices so that $0 \in Y_0^{(n)}$.

Let us denote by $P_n$ the collection of pieces of level $n$. The resulting
family $\mathcal{Y}_c := \{ P_0, P_1, \ldots \}$ of puzzle pieces of all
depths, has the following two properties:

{\renewcommand{\labelenumi}{{\bf P\arabic{enumi}}}
\renewcommand{\theenumi}{\labelenumi}
\begin{enumerate}
  \item \label{piece_relns} Any two puzzle pieces either are nested (with the
        piece of higher depth contained in the piece of lower depth), or have
        disjoint interiors.

  \item The image of any piece $Y_j^{(n)}$ $(n \geq 1)$ is a piece
        $Y_i^{(n-1)}$ of the previous depth $n-1$. The restricted map $f:{\rm
        int}\, Y_j^{(n)} \longrightarrow {\rm int}\, Y_i^{(n-1)}$ is a $2$ to
        $1$ branched covering or a conformal homeomorphism, depending on
        whether $j=0$ or not.
\end{enumerate}}

These properties characterize $\mathcal{Y}_c$ as a Markov family and endow the
puzzle partition with dynamical meaning.

Note that the collection of ray angles at depth $n$ consists of all
$n$-preimages of $\{ r_{t_1}, \ldots , r_{t_q} \}$ under angle doubling. The
union of all pieces of depth $n$ is the region enclosed by the equipotential
$e_{\left( 2^{2^{-n}} \right)}$. Note also that every piece $Y$ of depth $n$
is the $n^{\rm th}$ preimage of some piece of level 0. By further
iteration, $Y$ will map onto a region determined by the same rays as
$Y_0^{(0)}$ and a possibly larger equipotential. This provides a $1$ to $1$
correspondence between puzzle pieces and preimages of 0. The distinguished
point inside each piece is called the {\bf center} of the piece.

\subsection{Adjacency Graphs}\label{subsect:Graphs}
Given a set of puzzle pieces $P \subset P_n$, define the {\bf dual graph}
$\Gamma(P)$ as a formal graph whose set of vertices is $P$ and whose edges
join pairs of pieces that share an arc of external ray. Due to its finiteness,
it is always possible to produce an isomorphic model of $\Gamma(P)$ sitting in
the plane, without intersecting edges and such that it respects the natural
immersion of $\Gamma(P)$ in the plane.

\begin{defn}
  When $P = P_n$, we call $\Gamma_n:=\Gamma(P_n)$ the {\bf puzzle graph} of
  depth $n$. In this context, the vertices corresponding to the central piece
  $Y_0^{(n)}$ and the piece around the critical value $f_c(0)$ are denoted
  $\xi_n$ and $\eta_n$ respectively.
\end{defn}

\begin{defn}
  The vertices $\xi_n$ and $\eta_n$ determine two partial orders on the vertex
  set of $\Gamma_n$ as follows: If $a,b \in V(\Gamma_n)$, we write $a
  \succ_{\eta_n} b$ when every path from $a$ to $\eta_n$ passes through $b$.
  We write $a \succ_{\xi_n} b$ when every path from $a$ to $\xi_n$ passes
  through $b$ or through its symmetric image with respect to the origin.
\end{defn}

\begin{center}\begin{figure}[h]
  \includegraphics{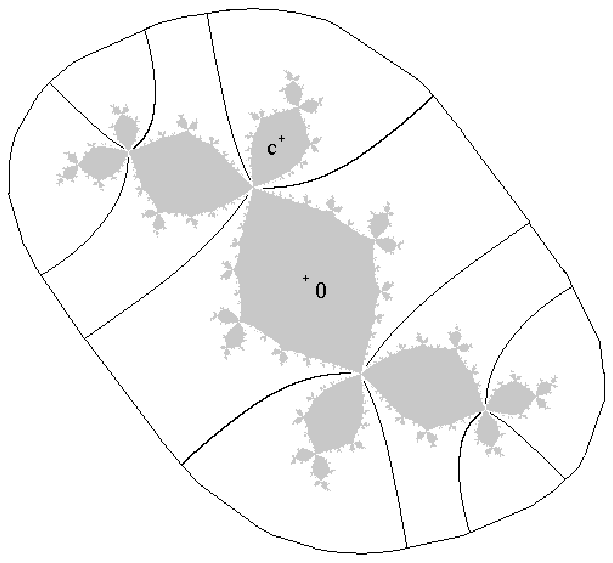}
  \includegraphics{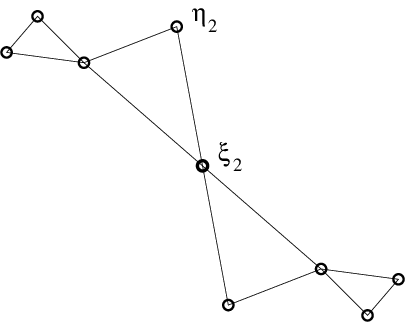}
  \caption[Puzzle of depth $2$ and its corresponding graph]
          {\label{fig:Puzzle_Graphs} \it Puzzle of depth $2$ and its
          corresponding graph. Splitting the graph at $\xi_2$ produces the
          graphs ${\rm Puzz}_2^-$ and ${\rm Puzz}_2^+$; both shaped like a
          bow tie and isomorphic to $\Gamma_1$. }
\end{figure}\end{center}

The following are natural consequences of the definitions; see Figure
\ref{fig:Puzzle_Graphs} for reference.

\begin{prop}\label{prop:Puzzle_Properties}
  The puzzle graphs of $f$ satisfy:

  {\renewcommand{\labelenumi}{{\bf G\arabic{enumi}}}
  \renewcommand{\theenumi}{\labelenumi}
  \begin{enumerate}
    \item \label{symmetry} $\Gamma_n$ has 2-fold central symmetry around
          $\xi_n$.

    \item \label{q-gon} $\Gamma_0$ is a $q$-gon whenever $c \in L_{p/q}$. For
          $n \geq 1$, $\Gamma_n$ consists of $2^n$ $q$-gons linked at their
          vertices in a tree-like structure; i.e. the only cycles on this
          graph are the $q$-gons themselves.

    \item \label{next_level} For $n \geq 1$, removing $\xi_n$ and its edges
          splits $\Gamma_n$ into $2$ disjoint (possibly disconnected)
          isomorphic graphs. Reattaching $\xi_n$ to each, and adding the
          corresponding edges defines the connected graphs ${\rm Puzz}_n^-$
          and ${\rm Puzz}_n^+$ (here, $\eta_n \in {\rm Puzz}_n^-$). Then
          $\Gamma_n = {\rm Puzz}_n^- \cup {\rm Puzz}_n^+$ and ${\rm Puzz}_n^-
          ,\, {\rm Puzz}_n^+$ are isomorphic to $\Gamma_{n-1}$ with $\mp
          \eta_n$ playing the role of $\xi_{n-1}$ in ${\rm Puzz}_n^{\pm}$.

    \item \label{maps} For $n \geq 1$ there are two natural maps:
          $f^*:\Gamma_n \longrightarrow \Gamma_{n-1}$ induced by $f$, and
          $\iota^*:\Gamma_n \longrightarrow \Gamma_{n-1}$ induced by the
          inclusion among pieces of consecutive depths. $f^*$ is $2$ to $1$
          except at $\xi_n$ and sends ${\rm Puzz}_n^{\pm}$ onto
          $\Gamma_{n-1}$. In turn, $\iota^*$ collapses the outermost $q$-gons
          into vertices.

    \item \label{order} The map $f^*:\big( \Gamma_n,\, \succ_{\xi_n} \big)
          \longrightarrow \big( \Gamma_{n-1},\, \succ_{\eta_{n-1}} \big)$
          respects order. That is, if $a \succ_{\xi_n} b$ then $f^*(a)
          \succ_{\eta_{n-1}} f^*(b)$.
  \end{enumerate}}
\end{prop}

\begin{defn}
  Let $\Gamma$ be a graph isomorphic to a subgraph of $\Gamma_n$ and $\Gamma'$
  a graph isomorphic to a subgraph of $\Gamma_{n-1}$. A map $E:\Gamma
  \longrightarrow \Gamma'$ that satisfies \ref{symmetry} and \ref{q-gon} will
  be called {\bf admissible} if it also respects order in the sense of
  \ref{order}.
\end{defn}

\begin{proof}[Proof of Proposition \ref{prop:Puzzle_Properties}.]
  Property \ref{symmetry} and the existence of $f^*$ and $\iota^*$ are
  immediate consequences of the structure of quadratic Julia sets. The
  configuration of $\Gamma_0$ is given by the rotation number around $\alpha$
  and then the tree-like structure of $\Gamma_n \, (n \geq 1)$ follows from
  \ref{next_level}.

  Consider a centrally symmetric simple curve $\gamma \subset Y_0^{(n)}$
  connecting two opposite points of the equipotential curve $e_{(2^{2^{-n}})}$
  that bounds $Y_0^{(n)}$. Then $\gamma$ splits the simply connected region
  $\bigcup_{Y \in P_n} Y$ in $2$ identical parts. Therefore, $\Gamma \setminus
  \xi_n$ is formed by $2$ disjoint graphs justifying the existence of
  ${\rm Puzz}_n^{\pm}$. However, $\bdry Y_0^{(n)}$ may contain several
  segments of $e_{(2^{2^{-n}})}$; so $\gamma$, and consequently
  ${\rm Puzz}_n^{\pm}$, {\it are not uniquely determined}. This ambiguity is
  not consequential; Lemmas \ref{lemma:New_Frame_Graph} and
  \ref{lemma:Embedding_of_F} describe the proper method of handling it.

  The fact that $f$ maps the central piece to a non-central one containing the
  critical value legitimizes the selection of ${\rm Puzz}_n^-$ as the unique
  graph containing $\eta_n$. By symmetry, every piece of $P_n$ except the
  central one has a symmetric partner and they both map in a $1$ to $1$
  fashion to the same piece of $P_{n-1}$. The isomorphisms in \ref{next_level}
  follow.

  If two pieces $A,B$ of depth $n$ share a boundary ray, their images will
  too. Moreover, letting $A',B'$ be the pieces of depth $n-1$ containing $A$
  and $B$, it is clear that $\bdry A'$ and $\bdry B'$ must share the {\it
  same} ray as $\bdry A$ and $\bdry B$. This shows that $f^*$ and $\iota^*$
  effectively preserve edges and are well defined graph maps. Clearly $f^*$ is
  $2$ to $1$, so to complete the proof of \ref{maps} it is only necessary to
  justify the collapsing property of $\iota^*$, and by Property
  \ref{next_level}, it is sufficient to consider the case $\iota^*:\Gamma_1
  \longrightarrow \Gamma_0$. Now, the non-critical piece $Y_j^{(0)}$ contains
  a unique piece $Y_j$ of $P_1$. However, the critical piece $Y_0^{(0)}$
  contains a total of $q$ different pieces of depth $1$: a smaller central
  piece $Y_0^{(1)}$ and $q-1$ lateral pieces $-Y_j$. The resulting graph,
  $\Gamma_1$, consists then of two $q$-gons joined at the vertex
  $\xi_1$. Under $\iota^*$, one of these $q$-gons collapses on the critical
  vertex $\xi_0$.

  To prove \ref{order}, let us construct the tree $\Gamma'_n$ with $2$ to $1$
  central symmetry by collapsing every $q$-gon into a single vertex. The
  orders $\succ_{\xi'_n},\, \succ_{\eta'_n}$ in $\Gamma'_n$ are induced by the
  orders in $\Gamma_n$. Then the corresponding map ${f^*}':\big( \Gamma'_n,\,
  \succ_{\xi_n} \big) \longrightarrow \big( \Gamma'_{n-1},\, \succ_{\eta_n}
  \big)$ is a $2$ to $1$ map on trees that takes each half of $\Gamma'_n$
  injectively into a sub-tree of $\Gamma'_{n-1}$ and respects order. Since
  vertices in a cycle are not ordered, $f^*$ respects order as well.
\end{proof}

\subsection{Parapuzzle}
While the puzzle encodes the combinatorial behavior of the critical orbit for
a specific map $f_c$, the {\it parapuzzle} dissects the parameter plane into
regions of parameters that share similar behaviors: In every wake of $M$ we
define a partition in pieces of increasing depths, with the property that all
parameters inside a given {\it parapiece} share the same critical orbit
pattern up to a specific depth.

\begin{defn}
  Consider a wake $W_{p/q}$ and let $n \geq 0$ be given. Call $W^n$ the wake
  $W_{p/q}$ truncated by the equipotential $e_{\left( 2^{2^{-n}} \right)}$ and
  consider the set of angles $\mathcal{P}_n(\frac pq) = \big\{ t \mid 2^nt \in
  \mathcal{P}\big( \frac pq \big) \big\}$ (compare Subsection
  \ref{subsect:External_Rays_n_Wakes}). The {\bf parapieces} of $W_{p/q}$ at
  depth $n$ are the closures of the components of $W^n \setminus \big\{r_t
  \mid t \in \mathcal{P}_n(\frac pq) \big\}$.
\end{defn}

\begin{note}
  Even though the critical value $f_c(0)$ is simply $c$, it will be convenient
  to write $c \in \Delta$ when $\Delta$ is a parapiece and $f_c(0) \in V$ when
  $V$ is a piece in the dynamical plane of $f_c$. In general, we will use the
  notation ${\rm OBJ}[c]$ to refer to dynamically defined objects ${\rm OBJ}$
  associated to a specific parameter $c$.
\end{note}

\begin{defn}
  When the boundary of a dynamical piece $A$ is described by the same
  equipotential and ray angles as those of a parapiece $B$, this relation is
  denoted by $\bdry A \circeq \bdry B$.
\end{defn}

\begin{defn}
  Let $c \in M$ be a parameter whose puzzle is defined up to depth $n$.
  Denote by ${\rm CV}_n[c] \in P_n[c]$ the piece of depth $n$ that contains
  the critical value: $f_c(0) \in{\rm CV}_n[c]$.
\end{defn}

A consequence of Formula \ref{eqn:Phi_M} is the well known fact that follows.
For a proof of the main statement, refer to \cite{DH_p-l} or
\cite{Roesch}. For a proof of the winding number property, refer to
\cite{Chirurgie} and Proposition 3.3 of \cite{L_parapuzzle}.

\begin{prop}\label{prop:Identify_Bdries}
  Let $\Delta$ be a parapiece of depth $n$ in some wake $W$. Then ${\rm
  CV}_n[c] \circeq \Delta$ for every $c \in \Delta$. The family $\big\{ c
  \mapsto {\rm CV}_n[c] \mid c \in \Delta \big\}$ is well defined and it
  determines a holomorphic motion of the critical value pieces. The
  holomorphic motion has $\big\{ c \mapsto f_c(0) \big\}$ as a section with
  winding number $1$.
\end{prop}

The result on winding number can be interpreted as loosely saying that, as $c$
goes once around $\bdry \Delta$, the critical value $f_c(0)$ goes once around
$\bdry{\rm CV}_n$. However, this description is not entirely accurate since
$\bdry{\rm CV}_n[c]$ changes with $c$.

Let us mention the following examples of combinatorial properties that depend
on the behavior of the first $n$ iterates of 0. The fact that these entities
remain unchanged for $c \in \Delta$ follows from Proposition
\ref{prop:Identify_Bdries} and will be useful in the next sections.

\begin{itemize}
  \item The isomorphism type of $\Gamma_n[c]$.
  \item The combinatorial boundary of every piece of depth $\leq n$.
  \item The location within $P_n[c]$ of the first $n$ iterates of the critical
        orbit. \\
\end{itemize}

From the general results of \cite{L_parapuzzle}, we can say more about the
geometric objects associated to the above examples.

\begin{prop}
  Each of the sets listed below moves holomorphically as $c$ varies in
  $\Delta$:

  \begin{itemize}
    \item The boundary of every piece of depth $\leq n$.
    \item The first $n$ iterates of the critical orbit.
    \item The collection of $j$-fold preimages of $\alpha$ and $\beta$ $(j
          \leq n)$.
  \end{itemize}
\end{prop}

\subsection{The principal nest}\label{subsect:Principal_Nest}
The principal nest is well defined for parameters $c$ that belong neither to
$\overline{\cardioid}$ nor to an immediate $M$ copy. The first condition means
that both fixed points are repelling (so the puzzle is defined), while the
second condition characterizes those polynomials that do not admit an {\it
immediate renormalization} as described below. We restrict further to
parameters $c$ such that the orbit of 0 is recurrent to ensure that the nest
is infinite. These necessary conditions will justify themselves as we describe
the nest.

\begin{center}\begin{figure}[h]
  \includegraphics{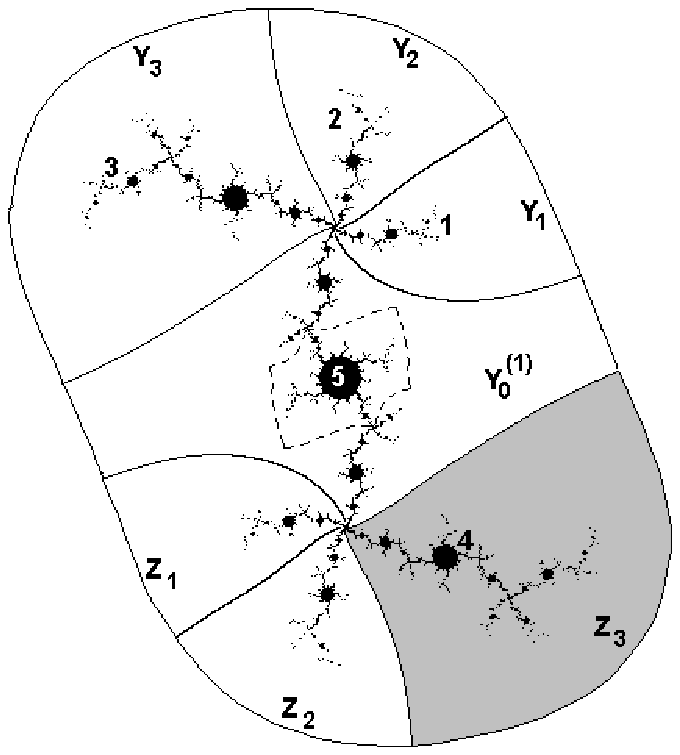}
  \caption[First stage of construction of the principal nest]
          {\label{fig:Nest_Of_Level_0} \it The puzzle $P_1(f_c)$ of depth $1$,
          where $c=(0.35926...)+i(0.64251...)$ is the center of the component
          of period $5$ in $L_{1/4}$. The first escape is $f_c^{\circ 4}(0) \in
          Z_3$ and the pull-back $V_0^0$ is shown in dotted lines. Note that
          $f^{\circ 5}(0) \in V_0^0$. This creates at once the piece $V_0^1
          \Subset V_0^0$ around the central component of $\mathbb{C} \setminus
          J_c$ ($V_0^1$ is not shown).}
\end{figure}\end{center}

In order to explain the construction of the principal nest, a more detailed
description of the puzzle partition at depth $1$ (use Figure
\ref{fig:Nest_Of_Level_0} for reference) is necessary. As a note of warning,
the pieces of depth $1$ will be renamed to reflect certain properties of
$P_1$. That is, we will override the use of the symbols $Y_j^{(1)}$. \\

The puzzle depth $P_1$ consists of $2q-1$ pieces of which $q-1$ are the
restriction to lower equipotential of the pieces $Y_1^{(0)}, Y_2^{(0)},
\ldots, Y_{q-1}^{(0)}$. Such pieces cluster around $\alpha$ and will be
denoted $Y_1, Y_2, \ldots, Y_{q-1}$. The restriction of $Y_0^{(0)}$ however,
is further divided into the union of the critical piece $Y_0^{(1)}$ and $q-1$
pieces $Z_1, Z_2, \ldots, Z_{q-1}$ which are symmetric to the corresponding
$Y_j$ and cluster around $-\alpha$. The indices are again determined by the
rotation number of $\alpha$ so that $f(Z_j)$ is opposite to $Y_j$ and
consequently $f(Z_j) = Y_{j+1}^{(0)}$.

Note that $f^{\circ q}(0) \in Y_0^{(0)}$, so we face two possibilities. It may
happen that $f^{\circ jq}(0) \in Y_0^{(1)}$ for all $j$, in which case it is
possible to find {\it thickenings} of $Y_0^{(1)}$ and $Y_0^{(0)}$, that yield
the {\bf immediate renormalization} $f^{\circ q}: Y_0^{(1)} \longrightarrow
Y_0^{(0)}$ described by Douady and Hubbard; or else, we can find the least $k$
for which the orbit of $0$ under $f^{\circ q}$ escapes from $Y_0^{(1)}$. We
will assume that this is the case, so $f^{\circ kq}(0) \in Z_\nu$ for some
$\nu$ and call $kq$ the {\bf first escape time.}

The initial nest piece $V_0^0$ is defined as the $(kq)$-fold pull-back of
$Z_{\nu}$ along the critical orbit; that is, the unique piece that satisfies
$0 \in V_0^0$ and $f^{\circ kq}(V_0^0) = Z_{\nu}$. In fact, $V_0^0$ can also
be defined as {\it the largest central piece that is compactly contained in}
$Y_0^{(1)}$: Notice that $Z_\nu \Subset Y_0^{(0)}$ so $V_0^0 \Subset
Y_0^{(1)}$; that is, $\big( {\rm int}\, Y_0^{(1)} \big) \setminus V_0^0$ is a
non-degenerate annulus.

The higher levels of the principal nest are defined inductively. Suppose that
the pieces $V_0^0, V_0^1, \ldots, V_0^n$ have been already constructed. If the
critical orbit never returns to $V_0^n$ then the nest is finite. Otherwise,
there is a first return time $\ell_n$ such that $f^{\circ \ell_n}(0) \in
V_0^n$; then we define $V_0^{n+1}$ as the {\it critical} piece that maps to
$V_0^n$ under $f^{\circ \ell_n}$.

\begin{prop} The principal nest $V_0^0 \Supset V_0^1 \Supset \ldots$ is a
  family of strictly nested pieces centered around 0.
\end{prop}

\begin{proof}
  $V_0^0$ is a piece of depth $kq$ (the first escape time). Since $V_0^1$ is a
  $f^{\circ \ell_1}$-pull-back of $V_0^0$, it is a piece of depth $kq+\ell_1$
  and, in general, $V_0^n$ will be a piece of depth
  $kq+\ell_1+\dots+\ell_n$. Since all pieces contain 0, Property
  \ref{piece_relns} implies that $V_0^j \supset V_0^{j+1}$.

  Recall that $V_0^0 \Subset Y_0^{(1)}$; thus, the $f^{\circ
  \ell_1}$-pull-backs of these $2$ pieces satisfy $V_0^1 \Subset X$ with $X$ a
  central piece of depth $1+\ell_1$. Now, $0 \nin Z_{\nu}$, so $f^{\circ
  kq}(0)$ requires further iteration to reach a central piece; i.e., $\ell_1 >
  kq$. By construction, $V_0^0$ is a central piece of depth $1+kq$, so
  Property \ref{piece_relns} implies $V_0^1 \Subset X \subset V_0^0$. An
  analogous argument yields the strict nesting property for the nest pieces of
  higher depth.
\end{proof}

\begin{defn}
  The {\bf principal annulus} $V_0^{n-1} \setminus V_0^n$ will be denoted
  $A_n$.
\end{defn}

It may happen that $\ell_{n+1} = \ell_n$; this means that not only does 0
return to $V_0^n$ under $f^{\circ \ell_n}$, but even deeper to $V_0^{n+1}$
without further iteration. In this case we say that the return is {\bf
central} and call a chain of consecutive central returns $\ell_n = \ell_{n+1}
= \ldots = \ell_{n+s}$ a {\bf cascade of central returns}. An infinite cascade
means that the sequence $\{ \ell_n \}$ is eventually constant, so $f^{\circ
\ell_n}(0) \in \bigcap_{j=n}^\infty V_0^j$. By definition, $f^{\circ
\ell_n}:V_0^{n+1} \longrightarrow V_0^n$ is a {\bf renormalization} of $f$;
that is, a $2$ to $1$ branched cover of $V_0^n$ such that the orbit of the
critical point is defined for all iterates.

The return to $V_0^n$, however, can be non-central. In fact, it is possible to
have several returns to $V_0^n$ before the critical orbit hits $V_0^{n+1}$ for
the first time. When a return is non-central, the description of the nest at
that level is completed by the introduction of the {\bf lateral} pieces $V_k^n
\in V_0^{n-1} \setminus V_0^n$. Let $\mathcal{O} \subset K$ denote the
critical orbit $\mathcal{O} = \left\{ f^{\circ j}(0)|j \geq 0 \right\}$ and
take a point $z \in \overline{\mathcal{O}} \cap V_0^{n-1}$ whose forward orbit
returns to $V_0^{n-1}$. If we call $r_{n-1}(z)$ the first return time of $z$
back to $V_0^{n-1}$, we can define $V^n(z)$ as the unique puzzle piece that
satisfies $z \in V^n(z)$ and $f^{\circ r_{n-1}(z)}\big(V^n(z)\big) =
V_0^{n-1}$. In particular, it is clear that $V^n(0)$ is just the same as
$V_0^n$ and that any $2$ pieces created by this process are disjoint or
equal.

\begin{defn}
  The collection of all pieces $V^n(z)$ where $z \in \overline{\mathcal{O}}
  \cap V_0^{n-1}$ that actually contain a point of $\mathcal{O}$ is denoted
  $\mathcal{V}^n$ and referred to as the {\bf level} $n$ of the nest.
\end{defn}

\begin{center}\begin{figure}[h]
  \includegraphics{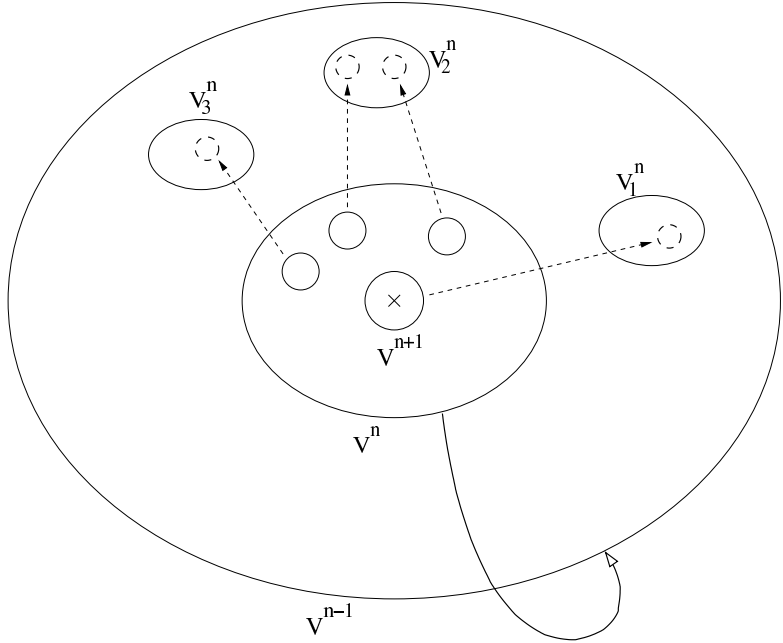}
  \caption[Relation between consecutive nest levels]{\label{fig:Sample_Nest}
           \it Relation between consecutive nest levels. The curved arrow
           represents the first return map $f^{\circ \ell_n}:V_0^n
           \longrightarrow V_0^{n-1}$ which is $2$ to $1$. The dotted arrows
           show a possible effect of this map on each nest piece of level
           $n+1$.  Each $V_j^{n+1}$ may require a different number of
           additional iterates to return to this level and map onto $V_0^n$.}
\end{figure}\end{center}

From this moment on, we will assume that the principal nest is infinite, with
$|\mathcal{V}^n| < \infty$ at every level, and that $f$ is non-renormalizable;
thus excluding the possibility of an infinite cascade of central returns. In
this situation we say that $f$ is {\bf combinatorially recurrent.} It follows
from \cite{L_attractor} and \cite{Martens} that $f$ acts minimally on the
postcritical set. In this situation, we can name the pieces $\mathcal{V}^n =
\{ V_0^n,V_1^n, \ldots ,V_{m_n}^n \}$ in such a way that the first visit of
the critical orbit to $V_i^n$ occurs before the first visit to $V_j^n$
whenever $i<j$. Obviously, the value of $r_{n-1}(z)$ is independent of $z \in
V_k^n$; thus we will denote it $r_{n,k}$.

\begin{defn}\label{defn:GQL_Maps}
  For every level of the nest, define the map:
  $$
    g_n:\bigcup_{k=1}^{m_n} V_k^n \longrightarrow V_0^{n-1},
  $$
  given on each $V_k^n$ by ${g_n}_{|_{V_k^n}} \equiv f^{\circ r_{n,k}}$.
\end{defn}

The map $g_n$ satisfies the properties of a {\it generalized quadratic-like}
({\bf gql}) {\it map,} i.e.:

{\renewcommand{\labelenumi}{{\bf
  g\hspace{-1pt}q\hspace{-1pt}l-\hspace{-1pt}\arabic{enumi}:}}
\renewcommand{\theenumi}{\labelenumi}
\begin{enumerate}
  \item $\bigcup_{\mathcal{V}^n} V_k^n \Subset V_0^{n-1}$ and all the pieces
        of $\mathcal{V}^n$ are pairwise disjoint.
  \item ${g_n}_{|_{V_k^n}} : V_k^n \longrightarrow V_0^{n-1}$ is a $2$ to $1$
        branched cover or a conformal homeomorphism depending on whether $k=0$
        or not.
\end{enumerate}}

Note that $g_n$ usually is the result of a different number of iterates of $f$
when restricted to different $V_k^n$. Since we often refer to the map $g_n$ as
acting on individual pieces, it is typographically convenient to introduce the
notation

\begin{defn}\label{defn:g_{n,k}}
  The map ${g_n}_{|_{V_k^n}} = f^{\circ r_{n,k}}$ will be denoted $g_{n,k}$.
\end{defn}

Thus, $g_{n,k}(V_k^n) = V_0^{n-1}$ is a $2$ to $1$ branched cover or a
homeomorphism depending on whether $k=0$ or not.

\subsection{Paranest}
The {\it paranest} is well defined around parameters $c$ outside the main
cardioid that are neither immediately renormalizable nor postcritically
finite.

\begin{defn}\label{defn:Paranest}
  If $c$ is a parameter such that $f_c$ has a well defined nest up to level
  $n$ (for $n \geq 0$), the {\bf paranest} piece $\Delta^n[c]$ is defined by
  the condition $\bdry \Delta^n[c] \circeq \bdry f_c(V_0^n)$; where $V_0^n$ is
  the central piece of level $n$ in the principal nest of $f_c$. By the
  Douady-Hubbard theory, $\Delta^n[c]$ is a well defined region.
\end{defn}

The definition of principal nest, together with Proposition
\ref{prop:Identify_Bdries} imply that when $c' \in \Delta^n[c]$, the principal
nests of $f_c$ and $f_{c'}$ are identical until the first return $g_n(0)$ to
$V_0^{n-1}$ (which creates $V_0^n$). In fact, the relevant pieces move
holomorphically as $c'$ varies and $\Delta^n[c]$ is the largest parameter
region over which the initial set of $\ell_n$ iterates of 0 (recall that $g_n
\equiv f^{\circ \ell_n}$) moves holomorphically without crossing piece
boundaries.

Following the presentation of \cite{L_parapuzzle}, the family
$\big\{g_n[c']:V_0^n[c'] \longrightarrow V_0^{n-1}[c'] \mid c' \in \Delta^n[c]
\big\}$ is a proper DH quadratic-like family with winding number 1. The last
property follows from Proposition \ref{prop:Identify_Bdries} since $g_n$ is
the first return to a critical piece at this level. \\

Since the central nest pieces are strictly nested, the above definition
implies that the pieces of the paranest are strictly nested as well. It
follows that $\big( {\rm int}\, \Delta^n \big) \setminus \Delta^{n-1}$ is a
non-degenerate annulus. One of our main concerns is to estimate its modulus
or, as it is sometimes called, the {\bf paramodulus.}

\section{Frame system}\label{sect:Frames}

Let $f_c$ have an infinite principal nest. We need a description of the
combinatorial structure around nest pieces in order to record their positions
relative to each other. In this Section we enhance the principal nest with the
addition of a {\it frame system.} The notion of frame, introduced in
\cite{1st_part}, provides the necessary language to locate the lateral nest
pieces and describe as a consequence, the behavior of the critical orbit. {\it
The idea is to split the central nest pieces in smaller regions by a procedure
that resembles the construction of the puzzle.}

\subsection{Frames}
Figure \ref{fig:Nest_Of_Level_0} provides a useful reference for the
construction of the initial frames $F_0,\,F_1$ and $F_2$. Some attention is
necessary at these levels to ensure that the properties of Proposition
\ref{prop:Frame_Properties} hold. Starting with level 3, frames are defined
recursively.

Consider the puzzle partition at depth $1$ and recall that $kq$ denotes the
first escape of the critical orbit to $Z_{\nu}$. The {\bf initial frame} $F_0$
is the collection of nest pieces $F_0 = \big\{Y_0^{(1)} \big\} \cup \big\{
\bigcup_{j=1}^q \{Z_j \} \big\}$; clearly, $\Gamma(F_0)$ is a $q$-gon. The
frame $F_1$ is the collection of $(f^{\circ kq})$-pull-backs of cells in $F_0$
along the orbit of 0.

From the definition, the central piece $V_0^0$ that maps $2$ to $1$ onto
$Z_{\nu} \in F_0$, is one of the cells of $F_1$. The pull-back of any other
cell $A \in F_0$ consists of two symmetrically opposite cells, each mapping
univalently onto $A$. We say that $F_1$ is a {\it well defined unimodal}
pull-back of $F_0$.

\begin{lemma}\label{lemma:Start_Frame}
  All the cells of $F_1$ are contained in $Y_0^{(1)}$.
\end{lemma}

\begin{proof}
  Since $kq > 1$, $f^{\circ kq}\big( Y_0^{(1)} \big)$ is an extension of
  $Y_0^{(0)}$ to a larger equipotential. Thus, $f^{\circ kq}\big( Y_0^{(1)}
  \big)$ contains all cells of $F_0$.
\end{proof}

Let $\lambda$ be the first return time of 0 to a cell of $F_1$. By Lemma
\ref{lemma:Start_Frame}, the collection $F_2$ of pull-backs of cells in $F_1$
along the $(f^{\circ \lambda})$-orbit of 0 is well defined and $2$ to $1$.

\begin{lemma}
  The frame $F_2$ satisfies:
  \begin{enumerate}
    \item All cells of $F_2$ are contained in $V_0^0$.
    \item $V_0^1$ is contained in the central cell of $F_2$.
  \end{enumerate}
\end{lemma}

\begin{proof}
  First note that $\lambda = kq + (q - \nu)$ is the first return of 0 to
  $Y_0^{(1)}$ after the first escape to $Z_{\nu}$. It follows that $kq <
  \lambda \leq \ell_0$, where the second inequality is true since $V_0^0 \in
  F_1$. Then the first return to $F_1$ occurs no later than the first return
  to $V_0^0$. By definition, $f^{\circ \lambda}(V_0^0)$ is just $Y_0^{(0)}$
  extended to a larger equipotential. Since all cells of $F_1$ are inside
  $Y_0^{(1)} \subset f^{\circ \lambda}(V_0^0)$, the first assertion follows.

  Now, $V_0^1$ is central. By the Markov properties of $\mathcal{Y}_c$, either
  $V_0^1$ is contained in the central cell $C$ of $F_2$ or vice versa.
  However, both $f^{\circ \ell_0}(V_0^1)$ and $f^{\circ \lambda}(C)$ belong to
  $F_1$. Since $\ell_0 \geq \lambda$, the first possibility is the one that
  holds. This proves property (2).
\end{proof}

\begin{figure}[h]
\begin{center}
  \includegraphics{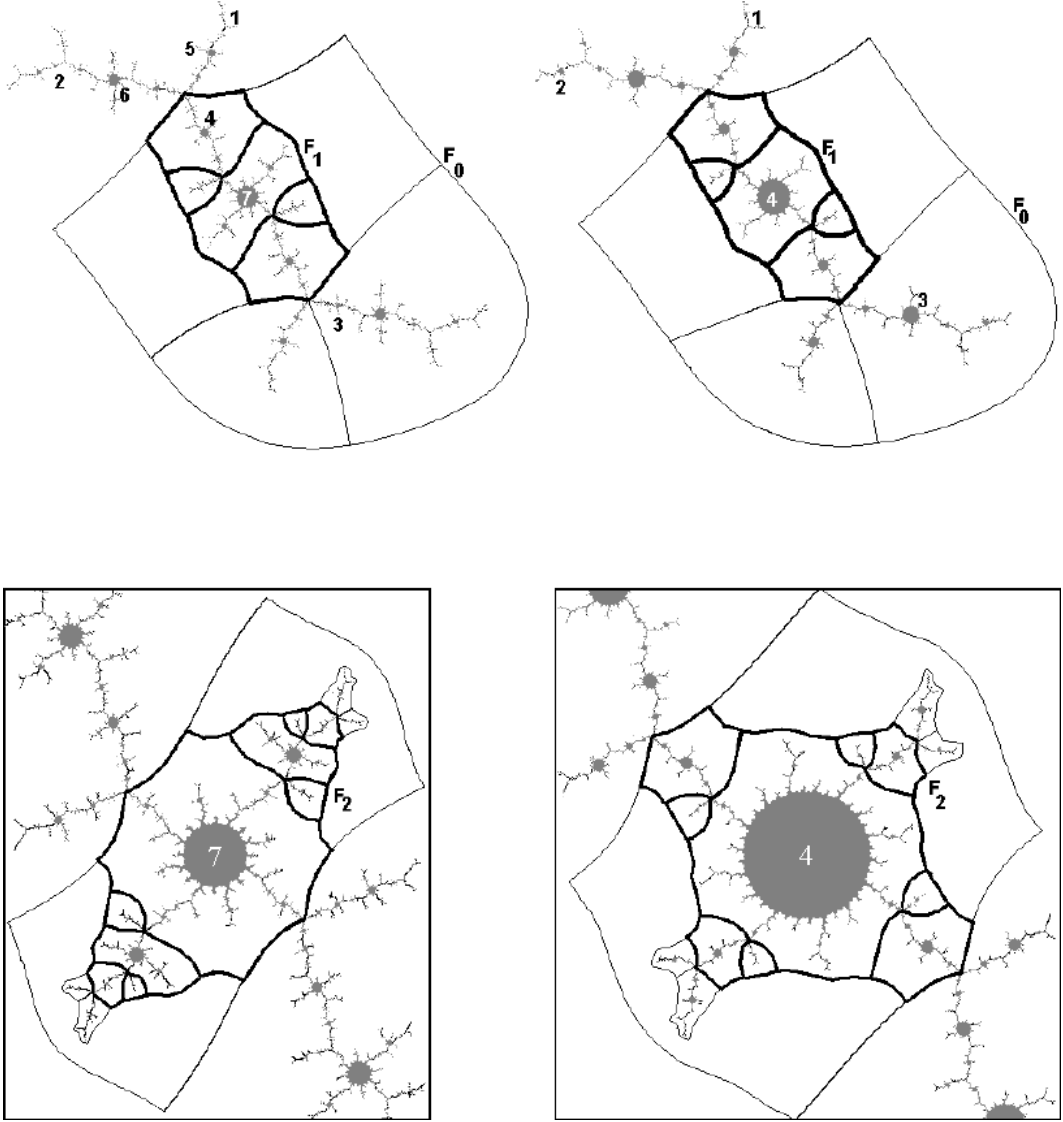}
  \caption[Comparison between two frames with similar combinatorics]
          {
           \label{fig:Two_Frames} \it Both of these parameters belong to the
           left antenna of $L_{1/3}$; they are centers of components of
           periods $7$ and $4$. Above we can see that the structures of the
           frames of levels 0 and $1$ coincide between the two
           examples. Still, the first return to $F_1$ falls in each case on a
           different cell, producing dissimilar frames of level 2. The
           pull-back of cells in $F_1$ produces the frame $F_2$, shown in
           heavy line on the second row.  }
\end{center}
\end{figure}

After introducing the first frames and linking them to the initial levels of
the nest, we can give the complete definition of the {\it frame system.}  The
driving idea of this discussion is that the internal structure of a frame
$F_{n+2}$, represented by the graph $\Gamma(F_{n+2})$, provides a
decomposition of $J_f \cap V_0^n$ that describes the combinatorial type of the
nest at level $n+1$.

\begin{defn}
  For $n \geq 0$ consider the first return $g_n(0) \in V_0^n$ and define
  $F_{n+3}$ as the collection of $g_n$-pull-backs of cells in $F_{n+2}$ along
  the critical orbit. The family $\mathcal{F}_c = \{F_0,F_1, \ldots \}$ is
  called a {\bf frame system} for the principal nest of $f_c$ and each piece
  of a frame is called a {\bf cell}.

  The dual graph $\Gamma(F_n)$ (see Subsection \ref{subsect:Graphs}) is called
  the {\bf frame graph.} As in the case of the puzzle graph, consider
  $\Gamma(F_n)$ with its natural embedding in the plane.
\end{defn}

Let us mention now some properties of frame systems (refer to
\cite{1st_part}).

\begin{prop} \label{prop:Frame_Properties} The frame system satisfies:
  \begin{enumerate}
    \item Frames exist at all levels.
    \item The central cell of $F_n$ contains the nest piece $V_0^{n-1}$.
    \item Each $F_n$ has 2-fold central symmetry around 0.
    \item \label{small_inside} Suppose there is a non-central return; then,
          eventually all nest pieces are compactly contained in cells of the
          corresponding frame.
  \end{enumerate}
\end{prop}

\subsection{Frame labels}\label{subsect:Frame_Lbls}
Our next objective is to introduce a labeling system for pieces of the
frame. This will allow us to describe the relative position of pieces of the
nest within a central piece of the previous level.  Unlike the case of
unimodal maps, where nest pieces are always located left or right of the
critical point, the possible labels for vertices of $\Gamma(F_n)$ will depend
on the combinatorics of the critical orbit. Only after determining the
labeling, it becomes possible to describe the location of nest pieces in a
systematic manner.

Observe that the structure of $F_{n+1}$ is determined by the structure of
$F_n$ and the location of $g_n(0)$. A graphic way of seeing this is as
follows. Say that the first return $g_{n-1}(0)$ to $V_0^{n-2}$ falls in a cell
$X \in F_n$. Let $L_n$ and $R_n$ be two copies of $\Gamma(F_n)$ with disjoint
embeddings in the plane. Now connect $L_n$ and $R_n$ with a curve $\gamma$
that does not intersect either graph. Suppose that one extreme of $\gamma$
lands at the vertex of $L_n$ that corresponds to $X$ and the other extreme
lands at the corresponding vertex of $R_n$ {\it approaching it from the same
access}.

\begin{figure}[h]
\begin{center}
  \includegraphics{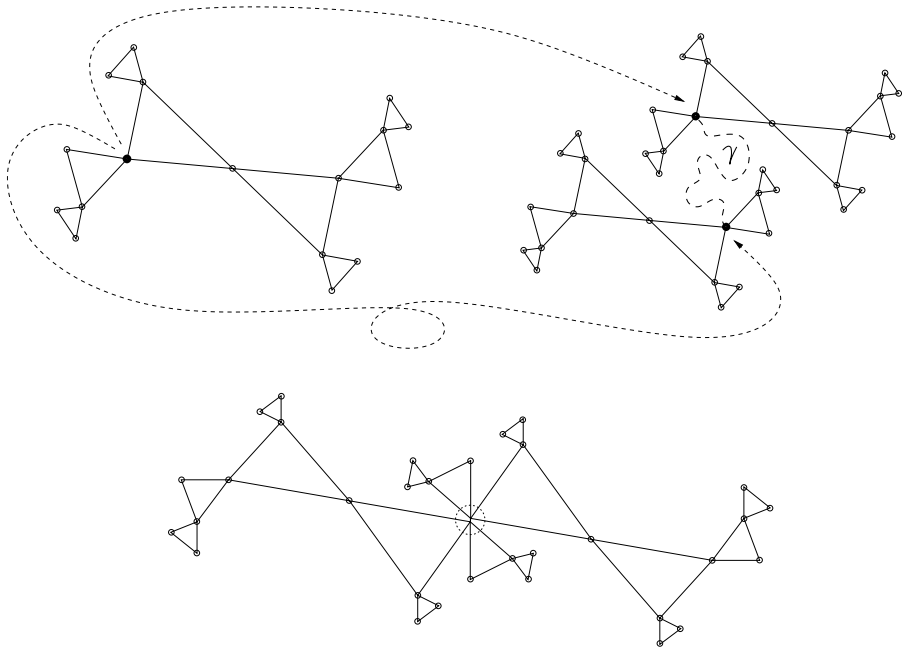}
  \caption[Constructing the next level graph.]{\label{fig:New_Frame} \it The
           curve $\gamma$ joins two copies of the same frame graph
           approaching the selected vertex from the same direction. The new
           frame graph is obtained after $\gamma$ is contracted to a point.}
\end{center}
\end{figure}

\begin{lemma}\label{lemma:New_Frame_Graph}
  If $\gamma$ is collapsed by a homotopy of the whole ensemble, the resulting
  graph is isomorphic to $\Gamma(F_{n+1})$.
\end{lemma}

\begin{lemma}\label{lemma:Embedding_of_F}
  The plane embedding of $\Gamma$ does not depend on the homotopy class of
  the curve $\gamma$ in lemma \ref{lemma:New_Frame_Graph}.
\end{lemma}

\begin{proof}
  Since we regard $\Gamma = \Gamma(F_n)$ as embedded in the sphere, the
  exterior of $\Gamma$ is simply connected, so there is a natural cyclic order
  of accesses to vertices (some vertices can be accessed from more than one
  direction). In this order, all accesses to $L_n$ are grouped together,
  followed by the accesses to $R_n$.
\end{proof}

A label at level $n$ will be a chain of $n+1$ symbols taken from the
alphabet $\{ {\sf Z}_0, {\sf Z}_1, \ldots, {\sf Z}_{(q-1)}, {\sf L}, {\sf R}
\}$. First, put the labels $\{ '{{\sf Z}_0}', '{{\sf Z}_1}', \ldots , '{{\sf
Z}_{(q-1)}}' \}$ on the cells of $F_0$, starting at the central piece
$Y_0^{(0)}$ and moving counterclockwise.

Let $\sigma_0$ be the label of the cell that holds the first return of 0 to
$F_0$ and, in general, let $\sigma_n$ denote the label of the cell in
$\Gamma(F_n)$ that holds the first return of 0. In order to label
$\Gamma(F_{n+1})$, assume that the number $q$ of pieces in $F_0$ is known, and
the {\it label sequence} $(q;\sigma_0, \ldots , \sigma_{n-1})$ that identify
the location of first returns of 0 to levels $0, \ldots , n-1$ of the nest.
In particular, all frames up to $\Gamma(F_n)$ have been successfully labeled.

Duplicate in $L_n$ the labels of $\Gamma(F_n)$, but concatenate an extra
$'{\sf L}'$ at the beginning. Do a similar labeling on $R_n$ by concatenating
an extra $'{\sf R}'$ to the duplicated labels. Note that the labels of the two
vertices corresponding to $X$ are $'{\sf L}'\sigma_n$ and $'{\sf
R}'\sigma_n$. The labels on $\Gamma(F_{n+1})$ will be the same as those in the
union of $L_n$ and $R_n$ except that we change the label of the identified
vertex, to become $'{{\sf Z}_0}'\sigma_n$. Clearly, $f$ induces a map
$f_*:\Gamma(F_{n+1}) \longrightarrow \Gamma_n$ for $n \geq 2$, that acts by
forgetting the leftmost symbol of each label.

It is important to mention that the resulting labeling of $\Gamma(F_n)$ {\bf
does} depend on the access to $\xi_n$ approached by $\gamma$. However, the
final unlabeled graphs are equivalent as embedded in the plane.

As was just mentioned, some vertices are accessible from $\infty$ in two or
more directions. These are precisely the vertices whose label contains the
symbol $'{{\sf Z}_0}'$ (for $n \geq 1$). Since such a vertex represents a
frame cell that maps (eventually) to a central frame cell, the tail of a label
with $'{{\sf Z}_0}'$ at position $j$ must be $\sigma_j$. On the other hand,
for every $j$ there must be labels with a $'{{\sf Z}_0}'$ in position $j$. It
follows that the set of labels of $\Gamma(F_n)$ and the sequence $(q;\sigma_0,
\ldots , \sigma_n)$ can be recovered from each other.

\subsection{Frame system and nest together}\label{sect:Frames_n_Nest}
The definition of frame system was conceived to satisfy the properties of
Proposition \ref{prop:Frame_Properties}. An extension of the argument used to
prove those properties shows that every piece $V_j^n$ of the nest is contained
in a frame cell of level $n+1$. Moreover, we would like to extend the
definition of frames so that each $V_j^n$ can be partitioned by a pull-back of
an adequate central frame. For this, recall first that $g_{n,j}(V_j^n) =
V_0^{n-1} \supset F_{n+1}$.

\begin{defn}
  The frame $F_{n,k}$ is the collection of pieces inside $V_k^{n-2}$ obtained
  by the $g_{n-2,k}$-pull-back of $F_{n-1}$. Elements of the frame $F_{n,k}$
  are called {\bf cells} and we will write $F_{n,0}$ instead of $F_n$, when
  there is a need to stress that a property holds in $F_{n,k}$ for every $k$.

  If a puzzle piece $A$ is contained in a cell $B \in F_{n,k}$, denote $B$ by
  $\Phi_{n,k}(A)$.
\end{defn}

We have described already how to label $F_n$. The other frames $F_{n,k} \, (k
\geq 1)$, mapping univalently onto $F_{n-1}$, have a natural labeling induced
from that of $F_{n-1}$ by the corresponding $g_{n-2,k}$-pull-back.

Let us describe now the itinerary of a piece $V_j^n$. Since $V_j^n \subset
V_0^{n-1}$, the map $g_{n-1}$ takes $V_j^n$ inside some piece
$V_{k_1(j)}^{n-1} \subset V_0^{n-2}$. Then, $g_{n-1,k_1(j)}$ takes
$g_{n-1}(V_j^n)$ inside a new piece $V_{k_2(j)}^{n-1}$ and so on, until the
composition of returns of level $n-1$
$$
  (g_{n-1,k_r(j)} \circ \ldots \circ g_{n-1,k_1(j)} \circ g_{n-1})_{|_{V_j^n}}
$$
is exactly $g_{n,j}:V_j^n \mapsto V_0^{n-1}$. Of course, $k_r$ is just 0, and
we will write it accordingly.

There is extra information that deems this description more accurate. For the
sake of typographical clarity, we will write $k_i$ instead of $k_i(j)$. For $i
\leq r$, let $\Phi_{n+1,k_i}$ be the cell in $F_{n+1,k_i} \subset
V_{k_i}^{n-1}$ that contains
$$
  g_{n-1,k_i} \circ \ldots \circ g_{n-1,k_1} \circ g_{n-1}(V_j^n)
$$
and denote by $\lambda_{n+1,k_i}$ the label of $\Phi_{n+1,k_i}$.

\begin{defn}
  The {\bf itinerary} of $V_j^n$ is the list of piece-label pairs:
  \begin{equation}
    \chi(V_j^n) = \left(
    [V_{k_1}^{n-1}; \lambda_{n+1,k_1}],
    [V_{k_2}^{n-1}; \lambda_{n+1,k_2}], \ldots,
    [V_{k_{r-1}}^{n-1}; \lambda_{n+1,k_{r-1}}],
    [V_0^{n-1}; \lambda_{n+1,0}] \right)
  \end{equation}
  up to the moment when $V_j^n$ maps onto $V_0^{n-1}$.
\end{defn}

Note first of all that the last label, $\lambda_{n+1,0}$, will start with
 $'{{\sf Z}_0}'$ due to the fact that $V_0^{n-1}$ is in the central cell of
 $F_n$. More importantly, the conditions
\begin{equation}\label{eqn:Admissibility}
  \begin{array}{lcll}
    V_{k_1}^{n-1} & \subset & g_{n-1}(\Phi_{n+1,0}) & \\
    V_{k_{i+1}}^{n-1} & \subset & g_{n-1,k_i}(\Phi_{n+1,k_i}) &
    2 \leq i <r
  \end{array}
\end{equation}
must hold since
$g_{n-1,k_{i-1}} \circ \ldots \circ g_{n-1,k_1} \circ g_{n-1}(V_j^n)
 \subset \Phi_{n+1,k_i}$ and
$g_{n-1,k_i} \circ \ldots \circ g_{n-1,k_1} \circ g_{n-1}(V_j^n)
 \subset V_{k_{i+1}}^{n-1}$.

\begin{defn}
  When the sequence of frame labellings is specified up to a given level $n$,
  the locations of the nest pieces and their (admissible) itineraries, we say
  that we have described the {\bf combinatorial type} of the map at level $n$.
\end{defn}

\section{$Q$-recurrency}\label{sect:Def_Q-rec}
Lyubich and Milnor established in \cite{LM_fibo} the uniqueness of the real
quadratic Fibonacci map $f_{c_{\rm fib}}$ and described in detail its
asymptotic geometry. The real parameter $c_{\rm fib} =-1.8705286321 \ldots$ is
determined by either of the following two equivalent conditions:

{\renewcommand{\labelenumi}{{\bf F\arabic{enumi}}}
\renewcommand{\theenumi}{\labelenumi}
\begin{enumerate}
  \item \label{closest-returns} The closest returns to 0 of the critical orbit
        occur exactly when the iterates are the Fibonacci numbers.
  \item \label{criss-cross} For $n \geq 2$, each level of the principal nest
        consists of the central piece $V_0^n$ and a unique lateral piece
        $V_1^n$. The first return map of previous level $g_{n-1}:V_0^{n-1}
        \longrightarrow V_0^{n-2}$ interchanges the central and lateral roles:
        $$g_{n-1}(V_0^n) \Subset V_1^{n-1}, g_{n-1}(V_1^n) = V_0^{n-1}.$$

        Additionally, the first returns to $Y_0^{(1)}$ and $V_0^0$ happen on
        the third and fifth iterates respectively.
\end{enumerate}}

The critical behavior of $f_{c_{\rm fib}}$ is the simplest among maps whose
nest has no central returns: Every level of the nest has a unique lateral
piece, so in a way, every first return comes as close as possible to being
central without actually being central.  This means that $f_{c_{\rm fib}}$ is
not renormalizable in the classical sense, although its combinatorics can be
described as an infinite cascade of {\it Fibonacci renormalizations} in the
space of {\bf gql} maps with one lateral piece.

The papers \cite{L_teichmuller} and \cite{W} analyze an unexpected feature of
the Fibonacci map. If the central pieces $V_0^n$ are rescaled to regions
$\Tilde{V}^n$ of fixed size, each $g_n$ induces a map $G_n:\Tilde{V}^n
\longrightarrow \Tilde{V}^{n-1}$. On increasing levels, the criss-cross
behavior that determines $c_{\rm fib}$ in condition \ref{criss-cross}
approximates with exponential accuracy the pattern of the critical orbit of
$P_{-1}(z) = z^2 - 1$ (i.e. $0 \mapsto -1 \mapsto 0 \mapsto -1 \mapsto
\ldots$). In fact, $G_n \longrightarrow P_{-1}$ locally uniformly in the $C^1$
norm.  Also, since ${\rm diam\,}\Tilde{V}^n \asymp 1$, it is shown that the
rescaled pieces converge in the Hausdorff metric to the filled Julia set of
$P_{-1}$.

In \cite{W}, Wenstrom translates this behavior to the Mandelbrot set and
obtains pieces of the paranest around $c_{\rm fib}$ that asymptotically
resemble $K_{-1}$; see Figure $1$ of \cite{W}. As consequences of this control
on shape, he computes the exact rate of linear growth of the principal moduli
and proves hairiness around the parameter $c_{\rm fib}$. \\

Let $c$ be the center of a prime hyperbolic component and $Q(z)$ its
associated polynomial. The critical orbit is periodic (of least period $m$)
and $Q^{\circ m}$ is the only renormalization of $Q$. An important consequence
of this, is that high enough depths of the puzzle of $Q$ will isolate in
individual pieces each point of the critical orbit $\mathcal{O}(Q) = \{0
\mapsto c \mapsto z_2 \mapsto \ldots \mapsto z_{m-1} \}$. Let us assume that
the fixed point $\alpha$ of $Q$ has combinatorial rotation number $\frac pq$.
In what follows we will save notation by restricting the use of ``$P_n$'' to
refer to the puzzle of $Q$ and ``$V_j^n$'', ``$F_n$'' for the nest and frames
of $Q$-recurrent maps.

Let us label $\Gamma(P_0)$, the graph of the puzzle of $Q$ at depth 0, with
symbols $'{{\sf Z}_0}'$ to $'{{\sf Z}_{q-1}}'$ starting at the critical point
piece and moving counterclockwise. Since $P_{n+1}$ is a $2$ to $1$ pull-back of
$P_n$, the graph $\Gamma(P_{n+1})$ consists of two copies of $\Gamma(P_n)$
identified at the critical value vertex and we can launch a labeling procedure
identical to the frame labeling of Subsection \ref{subsect:Frame_Lbls}. Note
that $\Gamma(P_n)$ is symmetric, but a canonical orientation can be specified
by dictating that the label on the critical value vertex begins with the
symbol $'{\sf L}'$. For $Q \in L_{p/q}$ the puzzle label sequence begins
$(q;\, '{{\sf Z}_p}', \ldots )$.

The above procedure creates a labeling of the puzzle of $Q$. Now consider any
map $f$ in the $\frac pq$-limb, with first escape time $q$ and such that
$f^{\circ q}(0) \in Z_p$. The map $f$ satisfies
\begin{itemize}
  \item The initial frame $F_0$ of $f$ consists of $q$ pieces and
        $\Gamma(F_0)$ is isomorphic to $\Gamma(P_0)$.
  \item The first return to $F_0$ is on the cell $Z_p$ which corresponds,
        under the above isomorphism, to the critical value piece of
        $P_0$. Therefore
  \item $\Gamma(F_1)$ is isomorphic to $\Gamma(P_1)$.
\end{itemize}
There is in fact, a full family $\Delta$ of parameters $c$ such that $f_c$
satisfies the above condition. Since the puzzle of $Q$ is created by
successive pull backs of the configuration $P_0$, the labeling of the puzzle
of $Q$ determines a weak admissible type in $\Delta$. Then Corollary 3.7 of
\cite{1st_part} guarantees the existence of parameters $c \in \Delta$ such
that the frame system of $f_c$ has the same structure as the puzzle of $Q$.

Observe that $F_n$ is symmetric so there are two choices for the homeomorphism
identifying $\Gamma(F_n)$ with $\Gamma(P_n)$. Once a frame orientation is
selected, we have an admissible label system.

\begin{defn} A critically recurrent polynomial $f_Q$ whose frame system has the
  same label sequence $(q; p, \sigma_1, \sigma_2, \ldots )$ as the puzzle of
  $Q$ is called {\bf $Q$-recurrent} if it satisfies the following additional
  condition. For any $n \geq 0$ and $2 \leq k \leq m-1$, the $k^{\rm th}$
  return to $V_0^n$ is the composition $( g_n \circ \ldots \circ g_{n+k-2}
  \circ g_{n+k-1} )$.
\end{defn}

\begin{note}
  There is an annoying offset between nest levels and frame levels.  Because
  of it, $V_0^n$ is contained in the central cell of $F_{n+1}$ and contains in
  turn the cells of $F_{n+2.}$ The notation suffers slightly when discussing
  return maps to several consecutive levels; hopefully this complication is
  balanced by the advantage of matching every frame level with the
  corresponding depth of the puzzle of $Q$.
\end{note}

\begin{prop}\label{prop:Q_Recurrent_Structure}
  For a $Q$-recurrent map every sufficiently high level $n$ of the nest has
  exactly $m$ pieces $V_0^n, V_1^n, \ldots, V_{m-1}^n$. For any $0\leq j \leq
  m-1$, $V_j^n$ is contained in the cell of $F_{n+1}$ corresponding to the
  piece in $P_n$ that contains $z_j$.
\end{prop}

\begin{figure}[h]
  \includegraphics{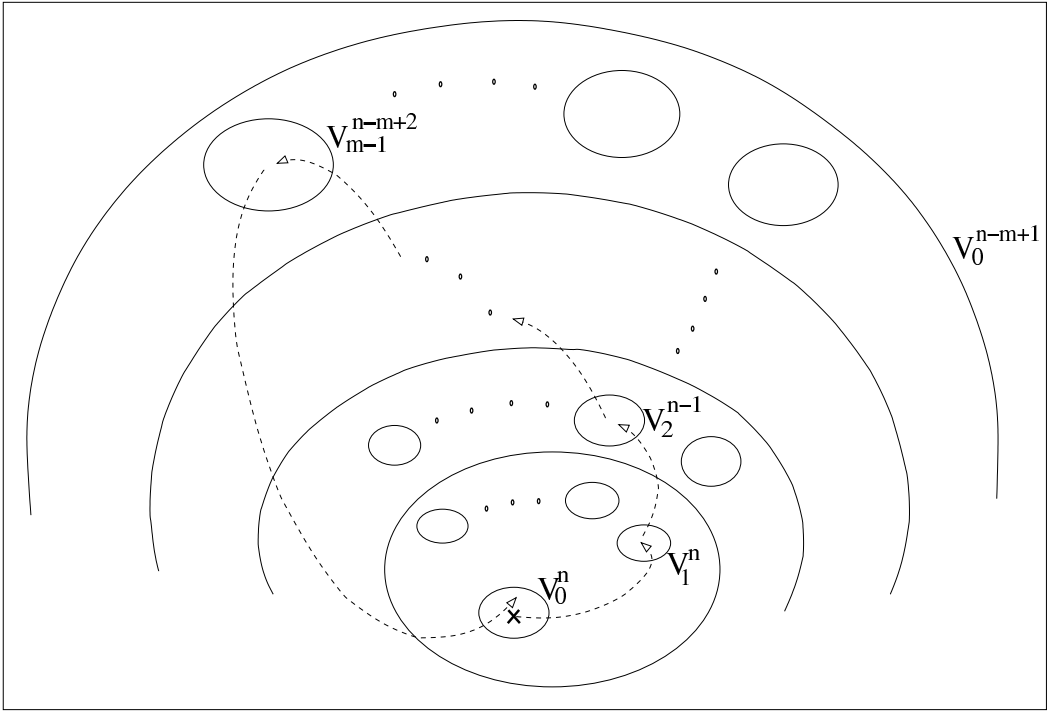}
  \caption[Nest of a $Q$-recurrent map]{\label{fig:Q_Recurrency} \it In a
           $Q$-recurrent map there are $m-1$ lateral pieces at each level,
           where $m$ is the period of 0 under $Q$. Under successive first
           return maps, the critical orbit jumps out to lower levels (as
           $V_i^{n-i+1}$ goes inside $V_{i+1}^{n-i}$) until at step $m-1$ it
           returns to the center. See also Figure \ref{fig:Detail_Of_Q-nest}.}
\end{figure}

\begin{proof}
  Choose $N$ big enough so that the puzzle $P_N$ isolates every point of
  $\mathcal{O}(Q)$ and let $n \geq N$. We will call $L_j^n$ the piece of $P_n$
  containing $z_j$.

  Consider the orbit of 0 under the composition $g_{n-2} \circ \ldots \circ
  g_{n+m-3}$. According to the label sequences, $g_{n+m-3}(0)$ falls in the
  cell of $F_{n+m-2}$ that corresponds to $L_1^{n+m-2}$. Next, $g_{n+m-4}\big(
  g_{n+m-3}(0) \big)$ falls in the cell of $F_{n+m-3}$ corresponding to
  $L_2^{n+m-3}$. Continue in this manner, with $g_{n+m-3-j} \circ \ldots \circ
  g_{n+m-3}(0)$ (where $0 \leq j \leq m-2$) falling in the cell of level
  $n+m-2+j$ that corresponds to $L_{j+1}^{n+m-2+j}$. At every step, jump out
  one nest level and create in the process (by adequate pull-backs) the nest
  pieces $V_1^{n+m-4}, V_2^{n+m-5}, \ldots, V_{m-1}^{n-2}$. Note that all
  these are lateral pieces since they are contained in a frame cell that is
  not central. In fact, $V_{j+1}^{n+m-4+j}$ is in the cell of $F_{n+m-2+j}$
  that corresponds to $L_j^{n+m-2+j}$; see Figure \ref{fig:Q_Recurrency}.

  The last map in this chain of compositions is $g_{n-2}$. It brings the
  critical orbit very nearly to the center, inside $V_0^{n+m-3}$. To see this,
  remember that the definition of $Q$-recurrency requires that the composition
  of maps $g_{n-2} \circ \ldots \circ g_{n+m-3}: V_0^{n+m-2} \longrightarrow
  V_0^{n+m-3}$ is the first return to $V_0^{n+m-3}$, i.e. $g_{n-2} \circ
  \ldots \circ g_{n+m-3}$ is exactly the map $g_{n+m-2}$. \\

  In summary, if a point of $\mathcal{O}(f)$ falls in a piece $V_{j}^n$ (for
  $j \leq m-1$), the next return falls inside $V_{j+1}^{n-1}$. If it falls on
  a piece $V_{m-1}^n$, the next return falls $m$ levels deeper, inside
  $V_0^{n+m-1}$ and is in fact, the first return to this piece. Repeating this
  procedure at $m-1$ consecutive levels creates various pieces of different
  levels. Among these, the $m-1$ lateral nest pieces of level $n$, each
  corresponding to a point $z_j$ ($1 \leq j \leq m-1$) of the critical orbit
  of $Q$.
\end{proof}

A first consequence of Proposition \ref{prop:Q_Recurrent_Structure} is the
fact that the itinerary of $V_j^n$ in the nest of level $n-1$ is
\begin{equation}\label{eqn:Q_Itineraries}
\begin{array}{lcl}
  \chi(V_j^n) & = &
   \big( [V_{j+1}^{n-1};\lambda_{n+1,j+1}], [V_0^{n-1};\lambda_{n+1,0}] \big)
   \, \text{ when } 0 \leq j <m-1 \\
  \chi(V_{m-1}^n) & = & \big( [V_0^{n-1};\lambda_{n+1,0}] \big).
\end{array}
\end{equation}

This hints to a similarity between the actions of $Q$ and $g_n$ that will be
made precise in the next Section, where we develop the asymptotic properties
of $Q$-recurrent maps and their principal nests.

\section{Asymptotics of $Q$-recurrency}\label{sect:Q_Recurrency}
This Section presents the geometric properties of $Q$-recurrent maps. Their
explicit relation to the combinatorics of the map $Q$ gives control over the
shapes of nest pieces and this will yield very precise estimates of the
analytic invariants of the nest.

For reference, let us state again the fundamental relation between levels of a
$Q$-recurrent nest:
\begin{eqnarray}
  g_{n+m} = g_n \circ \ldots \circ g_{n+m-1} \label{1st_Return_Composition} \\
  g_{n+m} \left( V^{n+m}, V^{n+m-1} \right) =
          \left( V^{n+m-1}, V^{n-1} \right) \label{Annuli_Returns}
\end{eqnarray}

We will also make repeated reference to the following result (see
\cite{L_nest}).

\begin{Lthm}
  Let $\kappa(n)$ count the levels of the principal nest up to level $n,$
  which are non-central. Then the moduli of the principal annuli grow
  linearly:
  $$
    \mu_{(n+2)} \geq B \cdot \kappa(n) + C.
  $$
  where the constant $B$ depends only on the initial modulus.
\end{Lthm}

\subsection{Complex Fibonacci maps}

Section \ref{sect:Def_Q-rec} begins with the definition of the real parameter
$c_{\rm fib}$. We will see that properties \ref{closest-returns} and
\ref{criss-cross} are shared by many complex maps.

\begin{defn}
  A complex polynomial map, or equivalently its corresponding parameter, is
  said to be {\bf complex Fibonacci} if the first return to each level of the
  nest happens exactly when the iterates are the Fibonacci numbers.
\end{defn}

\begin{note}
  The first return to a level can be viewed as a close return in a
  combinatorial sense; that is, a return to a small central piece. Since
  Lyubich's Theorem guarantees that central pieces decrease in size, the
  definition of complex Fibonacci parameter is equivalent to its metric
  analogue, property (\ref{closest-returns}) at the beginning of Section
  \ref{sect:Def_Q-rec}.
\end{note}

Our first result is a classification of the Fibonacci behavior in the complex
case.

\begin{defn}
  The set of all $(z^2-1)$-recurrent parameters is denoted $\Fibo$.
\end{defn}

\begin{thm}\label{thm:Fibonacci_Classification}
  A parameter $c$ is complex Fibonacci if and only if $c \in \Fibo$.
\end{thm}

\begin{proof}
  All $(z^2-1)$-recurrent maps have the same weak combinatorial type. As was
  pointed out in the note above, the first returns of high levels are just
  predetermined compositions of lower level ones. Thus, the number of iterates
  until the first return to a piece $V_0^n$ is independent of the parameter $c
  \in \Fibo$. Since the real parameter $c_{\rm fib} \in \Fibo$ is complex
  Fibonacci, the first direction of the assertion follows.

  To show the converse, it is only necessary to observe that the first return
  times in a Fibonacci nest are strictly increasing, so there are no central
  returns. Therefore the first return map $g_{n+1}$ must be the composition of
  at least two first return maps of the two preceding levels.  If the nest
  does not have $(z^2-1)$-recurrent type, there must be more than one lateral
  piece at some level $n$. Then the composition of maps generating $g_{n+1}$
  will actually contain more than two maps and the sequence $\{ \ell_n \}$ of
  first return times grows faster than the sequence generated by the recursion
  $\ell_{n+1} = \ell_n + \ell_{n-1}$. This contradicts the assumption that the
  map is complex Fibonacci.
\end{proof}

Although Yoccoz's Theorem on rigidity of non-renormalizable maps allows us to
characterize $\Fibo$ as a Cantor set, we must wait until next Section to show
that the relevant parapieces shrink exponentially fast, thus allowing us to
complete the description of the set $\Fibo$ as a Cantor set of Hausdorff
dimension 0 on which we can impose a natural dyadic decomposition.

\subsection{Shape}\label{sect:Shape}

We want to study the shape of nest pieces in the following sense.

\begin{defn}
  A sequence of compact sets $\{ C_j \subset \mathbb{C} \}$ is said to {\bf
  converge in shape} to a compact $K$ if there exist rescalings $\Tilde{C_j} =
  a_j \cdot C_j$ (with $a_j \in \mathbb{C}$) such that $\big\{ \Tilde{C_j}
  \big\} \longrightarrow K$ in the Hausdorff metric.
\end{defn}

The main Theorem of this Section is a vast extension of the result on the
shape of central pieces of $f_{c_{\rm fib}}$ found in \cite{L_teichmuller}. In
order to give the statement, some notation is needed.

Let $Q = Q(z)$ be the center of a prime hyperbolic component and $c_0$ a
$Q$-recurrent parameter. Recall that $f_{c_0}$ is described by a dyadic choice
of labels ($'{\sf L}'$ and $'{\sf R}'$) on every level. These frame
orientations determine the sequence of paranest pieces $\{ \Delta^n \}$ around
$c_0$. If $c \in \Delta^n$ is any nearby parameter, the combinatorics of $f_c$
are identical to those of $f_{c_0}$ including the orientations of the
homeomorphic frames, at all levels $j \leq n$. In particular, for any $c \in
\Delta^n$ we can find a (unique) point $s_j$ in $F_j$ corresponding to the
fixed point $\alpha$ of $Q$. In what follows, we omit from the notation the
fact that the objects described depend on $c$. Let $\alpha_j =
\frac{\alpha}{s_j}$ and define the complex rescalings $\Tilde{V}^j := \alpha_j
\cdot V_0^j$ of the central nest pieces of $f_c$, up to level $n$. Then, the
first return maps $g_j$ induce maps $G_j:\Tilde{V}^j \longrightarrow
\Tilde{V}^{j-1}$ on the rescaled pieces whose action on the rescaled frame
$\Tilde{F}_{j-2} := (\alpha_j \cdot F_{j-2}) \subset \Tilde{V}^j$ is isotopic
to the action of $Q$ on its own puzzle.

\begin{thm}\label{thm:Convergence_Of_G_n}
  Given $\varepsilon >0$, there is an $N$ such that for every parameter $c \in
  \Delta^n$ and level $n \geq N$, the maps $G_N,G_{N+1}, \ldots, G_n$ are all
  $\varepsilon$-close to $Q$ in the $C^1$ topology inside the ball of radius
  $\frac1{\varepsilon}$.
\end{thm}

\begin{corol}\label{corol:Convergence_In_Shape}
  The sequence of central nest pieces $\{ V_0^n \}$ of $f_{c_0}$ converges in
  shape to the filled Julia set $K_Q$.
\end{corol}

\begin{proof}[Proof of Corollary \ref{corol:Convergence_In_Shape}.]
  The point $\alpha \in K_Q$ is fixed under $G_n$ and is surrounded by
  $\Tilde{V}^n$. This rescaled nest piece also surrounds the critical point 0
  which attracts every point in $K_Q \setminus J_Q$. Now, $\Tilde{V}^n$ is the
  pull-back of $\Tilde{V}^{n-1}$ under $G_n$. By Theorem
  \ref{thm:Convergence_Of_G_n}, $G_n$ is a small perturbation of $Q$; since
  the rescaled pieces in the sequence $\{ \Tilde{V}^{n+1} , \Tilde{V}^{n+2} ,
  \ldots \}$ have bounded diameter, they become exponentially close to the
  regions in the sequence $\big\{ Q^{\circ -1}(\Tilde{V}^n), Q^{\circ
  -2}(\Tilde{V}^n), \ldots \big\}$ which converge to $K_Q$. This yields the
  result.
\end{proof}

\begin{figure}[h]
  \includegraphics{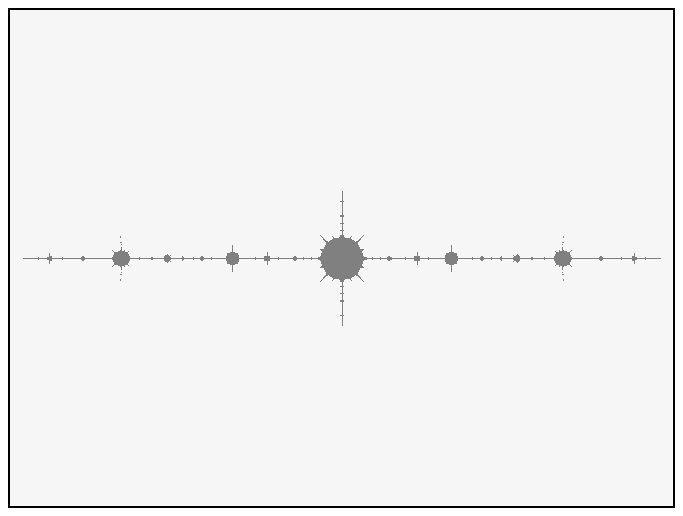} \hspace{6pt}
  \includegraphics{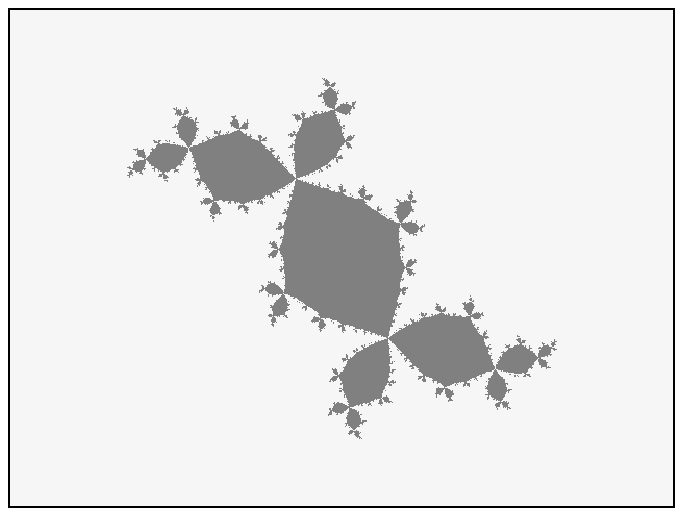}\\

  \includegraphics{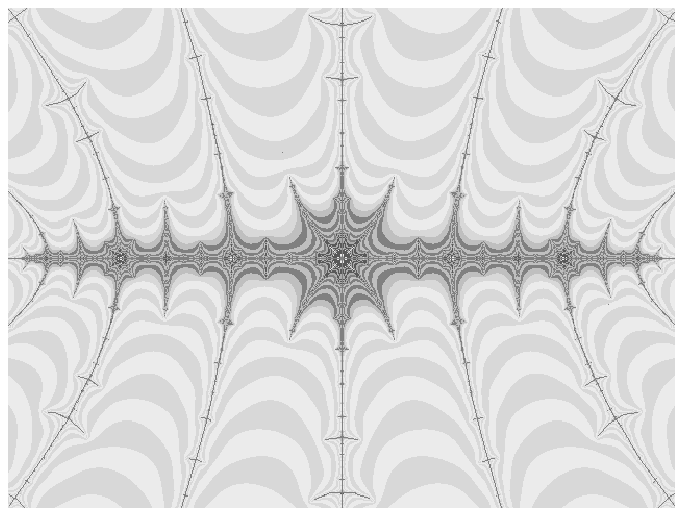} \hspace{6pt}
  \includegraphics{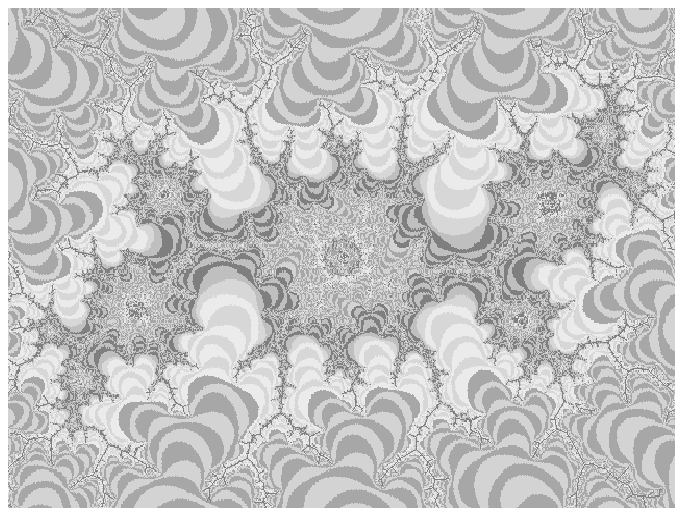}
  \caption[Deep levels of two $Q$-recurrent polynomials]
          {
           \label{fig:Period_3_Mimes} \it Consider the maps $Q_1:z \mapsto
           z^2-1.75487..$. (the ``airplane'') and $Q_2:z \mapsto
           z^2-(0.123...)+i(0.745...)$ (the ``rabbit''), as displayed in the
           first row. Both maps have critical orbits of period $3$. The
           pictures in the second row show close-ups near 0 of two other Julia
           sets. On the left, a $Q_1$-recurrent map
           ($c_1=-1.87449300898719..$.). On the right, a $Q_2$-recurrent map
           ($c_2=-0.023918090959967... + i0.984732550113053..$.). In each
           case, there is a central nest piece approximating the Julia set of
           the corresponding $Q_i$. }
\end{figure}

In particular, the central pieces of any quadratic complex Fibonacci map look
like $K_{-1}$, although each one may be tilted at a bizarre angle (recall that
the $\Tilde{V}^n$ are rescaled by a complex number). Other examples can be
seen in Figure \ref{fig:Period_3_Mimes}, showing puzzle pieces that
approximate the behavior of different periodic orbits of period $3$.

Notice that, since the frames are defined by the same sequence of pull-backs
as the central nest pieces, the result of Corollary
\ref{corol:Convergence_In_Shape} holds also for frames; i.e. the union of
cells in $F_n$ converges in shape to $K_Q$. \\

The proof of Theorem \ref{thm:Convergence_Of_G_n} depends on the convergence
of Thurston's map on an appropriate Teichm\"uller space (see the Appendix for
definitions). Let $\mathcal{O} \equiv \mathcal{O}(Q)$ and consider the surface
$S$ obtained by puncturing the plane at the critical orbit of $Q$; that is, $S
= \mathbb{C} \setminus \mathcal{O}$. Since deformations are considered only up
to an isotopy that leaves $\mathcal{O}$ invariant, the structure of a
puzzle-like construction does not change. Thus, when $h$ is a deformation in
the class of ${\rm id}$, the deformation $h\big( P(Q) \big)$ of the puzzle of
$Q$ can be isotoped back to the puzzle $P(Q)$ itself without changing its
configuration and without moving $\mathcal{O}$.

Thurston's map is best described via the alternate description of
$\mathcal{T}_S$ in terms of Beltrami differentials. First, normalize every
deformation $h$ by an affine change of coordinates $\varphi$ so that $\varphi
\circ h$ leaves $0,c \in \mathcal{O}$ fixed. The Beltrami coefficient $\mu =
\frac{\overline{\partial}h}{\partial h}\frac{d\overline{z}}{dz}$ determines a
conformal structure associated to $h$.

\begin{defn}
  The map $\tau_Q:\mathcal{T}_S \longrightarrow \mathcal{T}_S$ induced on
  equivalence classes of conformal structures by the pull-back $\mu \mapsto
  Q^*\mu$ is called the {\bf Thurston map} associated to $Q$.
\end{defn}

The action of $\tau_Q$ on a deformation class $h$ is easy to describe. The
class $\tau_Q\big( [h] \big)$ is represented by a deformation $\Tilde{h}$ such
that the map $Q_h = h \circ Q \circ \Tilde{h}^{-1}$ is analytic. Because of
conjugacy, $Q_h$ replicates the critical orbit behavior of $Q$ in a
neighborhood of $\Tilde{h}(\mathcal{O})$. In particular, one can specify a
puzzle-like structure around $\Tilde{h}(\mathcal{O})$ which pulls back
according to the same combinatorics as $Q$. Since $\mathcal{O}$ is finite, and
$Q^{\circ m}$ is not renormalizable, such a puzzle structure of high enough
depth will isolate all the elements of the critical orbit in individual
cells. We conclude that the isotopy class of $\Tilde{h}$ relative to punctures
consists of those $Q_h$-pull-backs of $h(\mathcal{O})$ that keep the puzzle
structure intact (however deformed).

\begin{proof}[Proof of Theorem \ref{thm:Convergence_Of_G_n}.]
  Let $X$ be any finite collection of simply connected compact subsets of
  $\mathbb{C}$. By a {\bf multicurve} $\Gamma$ around $X$ we mean a system of
  disjoint isotopy classes of simple closed curves in $\overline{\mathbb{C}}
  \setminus X$ such that each curve $\gamma_i \in \Gamma$ splits
  $\overline{\mathbb{C}}$ in two regions, each enclosing at least two elements
  of $X$ (i.e. $\gamma_i$ is non-peripheral).  If $f:\overline{\mathbb{C}}
  \setminus X \longrightarrow \overline{\mathbb{C}} \setminus X$ fixes every
  element of $X$, denote by $\Gamma_f^{-1}$ the multicurve consisting of the
  classes of $f$-preimages of elements $\gamma_i \in \Gamma$ that are not
  peripheral. The multicurve $\Gamma$ is said to be {\bf $f$-stable} if
  $\Gamma_f^{-1} \subset \Gamma$.

  Given a map $f:\overline{\mathbb{C}} \longrightarrow \overline{\mathbb{C}}$
  fixing the critical orbit $\mathcal{O}$ of $Q$ and an arbitrary $f$-stable
  multicurve $\Gamma$ around $\mathcal{O}$, we can construct the linear space
  $\mathbb{R}^{\Gamma}$ generated by the curves of $\Gamma$, and an induced
  linear map $\hat{f}_{\Gamma}:\mathbb{R}^{\Gamma} \longrightarrow
  \mathbb{R}^{\Gamma}$ given as follows. If $\gamma_i \in \Gamma$, let
  $\gamma_{i,j,k}$ denote the components of $f^{-1}(\gamma_i)$ that are in the
  class of $\gamma_j \in \Gamma_f^{-1}$. Then
  $$\hat{f}_{\Gamma}(\gamma_i) = \sum_{j,k} \frac1{d_{i,j,k}} \gamma_j.$$
  where $d_{i,j,k}$ denotes the degree of $f_{|_{\gamma_{i,j,k}}}:
  \gamma_{i,j,k} \longrightarrow \gamma_i$. \\

  An obstruction to the convergence of Thurston's map $\tau_f$ is represented
  by a $f$-stable multicurve around $\mathcal{O}$, for which
  $\hat{f}_{\Gamma}$ has an eigenvalue $\lambda \geq 1$. In our case, $Q$ is a
  polynomial so it represents the fixed point of its own Thurston map. In
  particular, there are no obstructions to the convergence of $\tau_Q$; see
  \cite{DH_thurston}.

  Now, since $Q$ belongs to a prime hyperbolic component of period $m$, the
  map $Q^{\circ m}$ is a renormalization conjugate to $z \mapsto z^2$. By
  hyperbolicity, the central puzzle pieces of $K_Q$ get arbitrarily close to
  the immediate basin of 0. In particular, there is a finite depth so that 0
  is the only point of $\mathcal{O}$ inside the central piece. By further
  iteration, the same will be true of any point in $\mathcal{O}$.

  Let us choose a level $k$ high enough so that the puzzle $P_{k-1}$ isolates
  all the points in the critical orbit of $Q$. Again, this is possible since
  $Q^{\circ m}$ is not renormalizable. Then, any $Q$-stable multicurve
  $\Gamma$ can be represented with curves that are constructed from segments
  of the arcs defining $P_{k-1}$. In this way, $\Gamma$ is described in terms
  of the structure of $P_{k'}$, for any level $k' \geq k-1$. Moreover,
  $\Gamma_Q^{-1}$ is a multicurve around $\mathcal{O}$ that can be described
  in terms of the combinatorial structure of $P_{k'+1}$. \\

  Now consider $f_c$ with $c \in \Delta^k$. Any $G_k$-stable multicurve
  $\Gamma'$ around the pieces $\Tilde{V}_j^{k+1}$ can be described with
  segments of curves in the boundary of the frame $\Tilde{F}_{k-1}$. Since
  $\Tilde{F}_{k-1}$ is isomorphic to $P_{k-1}$, there is a correspondence
  between $G_k$-stable multicurves around $\bigcup_j \left\{ \Tilde{V}_j^{k+1}
  \right\}$ and $Q$-stable multicurves around $\mathcal{O}$. This means that
  the only possible obstructions for $\tau_{G_k}$ must form inside one of
  those pieces; that is, a multicurve realizing such obstruction would
  intersect at least one of the pieces $\Tilde{V}_j^{k+1}$. Note that such
  multicurve cannot be represented by curves that are close to the boundary of
  $\Tilde{F}_{k-1}$.

  By \cite{L_nest}, the size of $\Tilde{V}_0^{k+2}$ with respect to
  $\Tilde{V}_0^{k+1}$ decreases exponentially as $k \rightarrow \infty$.
  Then, Koebe's Theorem implies that $G_k$ is exponentially close to being
  quadratic; that is, it can be decomposed as $G_k = D_k \circ Q_{h_k}$, where
  the maps $D_k$ become linear and the deformations $h_k$ are given by
  iteration of the Thurston map $\tau_Q$. Moreover, both $Q_{h_k}$ and $G_k$
  fix $\alpha$ and send 0 close to itself, so we can conclude that $D_k
  \longrightarrow {\rm id}$. It follows that $G_k$ rapidly approaches
  $Q_{h_k}$.

  Select any $Q$-stable multicurve $\Gamma'$ around $\mathcal{O}$. If there is
  a level $k$ such that $\Gamma'$ does not intersect any of the pieces
  $\Tilde{V}_j^{k+1}$, then $\Gamma'$ can be pushed to the boundary of
  $\Tilde{F}_{k-1}$ to represent a $G_k$-stable multicurve around the pieces
  $\Tilde{V}_j^{k+1}$. Since $\Gamma'$ is not an obstruction for $Q$, we
  deduce that, outside the $\Tilde{V}_j^{k+1}$, the map $G_k$ is isotopic to
  $Q$. However, the only possible Thurston obstructions are restricted to
  extremely small regions, then the distortion of $h_k$ goes to 0 and the maps
  $G_k$ converge to $Q$ exponentially fast in a neighborhood of $K_Q \setminus
  \mathcal{O}$. The Koebe space between $V_0^{k+1}$ and $V_0^k$ increases
  without bound, so we can claim convergence of the maps $G_k$ in arbitrarily
  big neighborhoods of $K_Q$.
\end{proof}

Theorem \ref{thm:Convergence_Of_G_n} has broad implications since it provides
excellent control of the shapes of nest pieces. In the next Subsection, we use
our knowledge on the shape of the central pieces to compute the rate of growth
of the principal moduli.

\subsection{Growth of annuli}
Here we study the moduli of the principal annuli in $Q$-recurrent maps. For
this family, we can state a more precise version of Lyubich's Theorem on the
linear growth of moduli. The key ingredient in our proof is Theorem
\ref{thm:Convergence_Of_G_n} giving control over the shape of pieces, together
with the extended Gr\"otzsch inequality as stated in the Appendix.

As a preparation for Theorem \ref{thm:Moduli_Growth}, we compute first the
capacities of $K_Q$ with respect to 0 and $\infty$.

\begin{lemma}
  Let $Q = Q(z)$ be the center of a hyperbolic component with critical period
  $m$. Then ${\rm cap}_{\infty} (K_Q) = 0$ and
  $$
    {\rm cap}_0 (K_Q) =
    -(m-1)\ln 2 - \sum_{j=1}^{m-1} \ln \left|Q^{\circ j}(0) \right|.
  $$
\end{lemma}

\begin{proof}
  $K_Q$ is connected, so the B\"ottcher coordinate $\varphi:
  \overline{\mathbb{C}} \setminus K_Q \longrightarrow \overline{\mathbb{C}}
  \setminus \mathbb{D}$ sending 0 to $\infty$ is the Riemann mapping with
  derivative 1, so the first equality is obvious.

  The capacity of $K_Q$ at 0 is simply ${\rm cap}_0 (U)$, where $U$ is the
  immediate basin of attraction of 0. Consider the iterated polynomial
  $Q^{\circ m}:U \longrightarrow U$. It is a $2$ to $1$ map of a simply
  connected domain with fixed critical point. Therefore, there is a map
  $\psi:\mathbb{D} \longrightarrow U$ such that
  \begin{equation}\label{eqn:Conjugate_Q^m}
    \psi(z^2) = Q^{\circ m} \circ \psi(z)
  \end{equation}
  and it is clear that
  ${\rm cap}_0 (K_Q) = {\rm cap}_0 (U) = \ln |\psi'(0)|$.

  Equation \ref{eqn:Conjugate_Q^m} shows that $\psi'(0)$ is the inverse of the
  quadratic coefficient in the series expansion of $Q^{\circ m}(z)$.  Since
  the constant term of $Q^{\circ j}(z)$ is just $Q^{\circ j}(0)$, it is easy
  to find recursively that $\frac1{\psi'(0)} = \prod_{j=1}^{m-1}2Q^{\circ
  j}(0)$ and thus,
  $$
    {\rm cap}_0(K_Q) = \ln\psi'(0) =
    -(m-1)\ln 2 - \sum_{j=1}^{m-1} \ln \left|Q^{\circ j}(0) \right|.
  $$
\end{proof}

Recall that capacity and modulus are invariants that vary continuously with
respect to the Carath\'eodory topology. Thus, given a sequence of topological
disks around 0 converging in the Hausdorff topology to a set with pinched
points, the sequence of capacities will converge to the capacity of the
component of the limit set that contains 0. Similarly, for a sequence of
annuli with adequate convergence, the limit of moduli detects only the modulus
of the limit component that contains the limit closed geodesic.

\begin{thm}\label{thm:Moduli_Growth}
  For any parameter $c \in \Fibo$ the principal moduli grow linearly at the
  rate
  $$\lim_{n \rightarrow \infty} \frac{\mu_n}n = \frac{\ln 2}3.$$

  If $Q = Q(z)$ is the center of a prime hyperbolic component with critical
  period $m \geq 3$, the rate of growth is exponential
  $$\lim_{n \rightarrow \infty} \frac{\mu_n}{\mu_{n+1}} = \kappa_m,$$ depends
  only on the period $m$ of $Q$ ad=nd satisfies $\kappa_m \nearrow \frac32$ as
  $m$ increases.
\end{thm}

\begin{proof}
  Fix a level $N$ large enough so that the shape of rescaled nest pieces is
  already close to the shape of $K_Q$. In particular, ${\rm cap}_0 \big(
  \Tilde{V}^n \big) \sim {\rm cap}_0 \big( (K_Q)_0 \big)$ for all $n \geq N$,
  where $(K_Q)_0$ is the Fatou component of $K_Q$ containing 0.

  We also require the lateral pieces are small enough to sit in the center of
  their (almost pinched) regions, far away from the boundary. This is possible
  since Lyubich's Theorem on linear growth forces shrinking and Theorem
  \ref{thm:Convergence_Of_G_n} locates the nest pieces in positions that
  resemble the critical orbit of $Q$.

  Theorem \ref{thm:Convergence_Of_G_n} also gives the recursion formula
  $$g_{n+m} = g_n \circ \ldots \circ g_{n+m-2} \circ g_{n+m-1}.$$

  On consecutive levels the first returns of a central piece $V_0^{n+m+1}$
  fall inside the pieces $V_1^{n+m}, V_2^{n+m-1},\ldots, V_{m-1}^{n+2}$ and
  $V_0^{n+m}$. From this, we obtain the annuli relation
  $$ g_{n+m}^{-1}\left( V_0^n \setminus V_0^{n+m} \right) = \left( V_0^{n+m}
    \setminus V_0^{n+m+1} \right).
  $$

  In order to estimate the modulus of $\left( V_0^{n+m} \setminus V_0^{n+m+1}
  \right)$, let us split the return map ${g_{n+m}}_{|_{ \left( V_0^{n+m}
  \setminus V_0^{n+m+1} \right) }}$ in the above mentioned composition of
  first returns. First, $g_{n+m-1}$ is $2$ to $1$ on the annulus $\left(
  V_0^{n+m} \setminus V_0^{n+m+1} \right)$. Note that the image of
  $V_0^{n+m+1}$ is deep inside $V_1^{n+m}$; in fact, these two pieces are
  separated by a nested sequence of preimages of the central pieces
  $V_0^{n+1},\ldots, V_0^{n+m-1}$. Due to the pinching of pieces near
  repelling points, most of the modulus of $g_{n+m-1}\left( V_0^{n+m}
  \setminus V_0^{n+m+1} \right) \subset V_0^{n+m-1}$ is concentrated in a
  region of $V_0^{n+m-1}$ that resembles the immediate basin of the critical
  value of $Q$. On this region, $g_{n+m-2}$ is injective.

  The remaining returns $g_{n+m-3},\ldots, g_n$ behave in a similar manner,
  essentially preserving the modulus on regions around non-central pieces. We
  can conclude that
  $$
    {\rm mod}\left( V_0^{n+m} \setminus V_0^{n+m+1} \right) \asymp
    \frac12 {\rm mod}\left( V_0^n \setminus V_0^{n+m} \right).
  $$

  The right hand side can be estimated by applying the Extended Gr\"otzsch
  Inequality to the annulus
  $\left( V_0^n \setminus V_0^{n+m} \right) =
   \left( V_0^n \setminus V_0^{n+1} \right) \cup
   \ldots \cup \left( V_0^{n+m-1} \setminus V_0^{n+m} \right)$, we obtain
  $$
    {\rm mod}\left( V_0^{n+m} \setminus V_0^{n+m+1} \right) \asymp
    \frac12 \sum_{j=0}^{m-1} {\rm mod}(V_0^{j} \setminus V_0^{j+1}) +
    \varepsilon_Q
  $$ where $\varepsilon_Q = (m-1) \left| {\rm cap}_0(K_Q) \right| = (m-1)
  \left| (m-1)\ln 2 + \sum_{j=1}^{m-1} \ln \left|Q^{\circ j}(0) \right|
  \right|$ \\

  The recursive formula $x_{n+m} = \frac12 (x_{n+m-1} + x_{n+m-2} + \ldots
  x_n) + \varepsilon_Q$ has an asymptotic behavior that is ruled by the
  largest real root of its characteristic polynomial
  \begin{equation}\label{eqn:Char_Pol}
    p(z) = z^m - \frac12 \Big( z^{m-1} + z^{m-2} + \ldots + z + 1 \Big)
  \end{equation}
  or, equivalently
  {\addtocounter{equation}{-1}
  \renewcommand{\theequation}{\arabic{equation}'}
  \begin{equation}\label{eqn:Char_Pol_prime}
    p(z) = z^m + \frac12 \frac{1-z^m}{z-1}.
  \end{equation}}

  When $m = 2$, the largest root of $z^2-\frac12(z+1)$ is 1. Consequently, the
  growth of the moduli is dominated by a linear term $\mu_n \sim An+B$. The
  only map with critical orbit of period $m=2$ is $Q(z)=z^2-1$, for which
  $\varepsilon_Q =\ln2$. The recursive relation $\mu_n \asymp
  \frac{\mu_{n-1}}2 + \frac{\mu_{n-2}}2 + \ln 2$ gives
  $$A = \lim_{n \rightarrow \infty} \frac{\mu_n}n = \frac{\ln 2}3.$$

  The analysis is different when $m \geq 3$. If $z \geq \frac32$, the second
  term of (\ref{eqn:Char_Pol_prime}) is $\frac12 \frac{1-z^m}{z-1} \leq
  1-z^m$, so the characteristic polynomial $p$ satisfies $p(z) \geq 1$. On the
  other hand, it follows from (\ref{eqn:Char_Pol}) that $p(1)<0$ (since $m
  \geq 3$). Thus the largest root $\kappa$ of $p(z)$ is in the interval
  $(1,\frac32)$ and the exponential growth of $\mu_n$ follows.

  Clearly, $\kappa$ does not depend asymptotically on the value of
  $\varepsilon_Q$. Nevertheless, $\kappa$ does vary with the period $m$ of
  $Q$. In fact, the sequence $\big\{ \kappa(m) \big\}$ converges to $\frac32$
  as $m$ increases. To see this, note from (\ref{eqn:Char_Pol_prime}) and
  $\kappa > 1$ that $p(\kappa)=0$ implies $2\kappa^{m+1}-3\kappa^m+1=0$. The
  polynomial $\hat{p}(z) = 2z^{m+1}-3z^m+1$ satisfies $\hat{p}(1)=0,\,
  \hat{p}(\frac32)=1$ and $\hat{p}'(1) <0$ for all $m \geq 3$. Since
  $\hat{p}'(z) = 2(m+1)z^m-3mz^{m-1}$, the interval $(1,\frac32)$ contains one
  critical point $z = \frac{3m}{2(m+1)}$ which squeezes $\kappa(m)$ towards
  $\frac32$.
\end{proof}

\begin{note}
  One should contrast the above result with \cite{Torrential}. There, the
  authors show that for almost every non-hyperbolic real parameter, the
  principal moduli grow at least as fast as a tower of exponentials. The
  ``slower'' growth displayed by $Q$-recurrent polynomials has immediate
  geometric consequences.
\end{note}

\begin{defn}
  Say that a compact set $K$ is {\bf hairy} at a point $c \in K$ if there is a
  sequence $\{ \varepsilon_1, \varepsilon_2, \ldots \}$ converging to 0, such
  that $\frac1{\varepsilon_j} \cdot (K-c) \cap \overline{\mathbb{D}}$ becomes
  dense in $\overline{\mathbb{D}}$.

  If $K$ is hairy at $c$ for any sequence of scaling factors $\{ \varepsilon_j
  \}$, we say that it satisfies {\bf hairiness at arbitrary scales}.
\end{defn}

 By an observation of Rivera-Letelier (\cite{Juan}), the construction of
 \cite{W} can be extended to prove hairiness of $M$ at any critically
 recurrent non-renormalizable parameter. The idea is as follows:

Since $K_c$ is connected, it contains a path from 0 to $\beta$. This crosses
every principal annuli from one boundary component to the other. Choosing a
high enough level, the annulus $A_n$ can be rescaled to constant diameter
containing a {\it hair} that connects the outer boundary with a small
neighborhood of 0. The pull-backs by consecutive first return maps duplicate
the number of hairs inside deeper annuli and this collection of hairs is
equidistributed around 0 (control of geometry). The hairiness of $K_c$ is then
translated to the parameter plane to obtain the result.

Rivera-Letelier has announced a proof that the real quadratic Fibonacci
polynomial displays hairiness at arbitrary scales. The proof relies in an
essential way on the linear growth of moduli, so it holds true for any
parameter in $\Fibo$. Other $Q$-recurrent polynomials miss this sharper
property in account of the exponential growth of their principal moduli. It
should be observed that this same property creates a somewhat embarrassing
difficulty; since the moduli grow so fast, computer generated pictures fail to
exhibit a convincingly hairy picture. In order to do so, it would be necessary
to reach deep levels of the nest that may be out of the range of resolution of
the software used.

\subsection{Meta-Chebyshev}

Starting with the chain ${\rm L} \check{\rm R} \check{\rm R} {\rm L}
\check{\rm L}$, construct an infinite sequence by the following iterative
procedure. At each step, concatenate a second copy of the current chain on
which the second to last marked symbol is substituted by its opposite. The
result is
$$
  \Theta:\, {\rm L} \check{\rm R} \check{\rm R} {\rm L} \check{\rm L} {\rm
  LRLL} \check{\rm L} {\rm LRRLRLRLL} \check{\rm L} {\rm LRRLLLRLLRLRRLRLRLL}
  \check{\rm L} \cdots
$$

In order to verify that this is an admissible kneading sequence, we have to
describe the sequence $\epsilon_1 \epsilon_2 \ldots$ of accumulated
orientation reversals in $\Theta$; that is, $\varepsilon_j$ is $+$ or $-$
depending on whether the number of L's up to position $j$ is even or odd. Then
we only have to prove that for any $m$, the least $i$ such that
$\epsilon_{m+i} \neq \epsilon_m \cdot \epsilon_i$ satisfies $\epsilon_i =\,
-$. The sequence of $\epsilon_j$ begins:
$$- \check- \check- + \check- + + - + \check- + + + - - + + - +
  \check- + + + - + - - + - - + + + - - + + - + \check- \cdots
$$ and the rule to construct it is as follows: Start with the chain $- \check-
\check- + \check-$. At each step make a second unchecked copy of the current
chain and invert every symbol that appears to the left of the second to last
check; then concatenate this copy to the right and put a check on the last
symbol.

This sequence starts with $- \check- \check- + \check-$ and every mark will be
on a $-$ symbol. It follows that there cannot be more that three $+$ in a row
so the admissibility condition is satisfied. Thus, the kneading sequence
$\Theta$ can be realized by a real polynomial. In fact, by $P_{\rm MCheb}(z)
:= z^2 - 1.87450961730020085 \ldots$

This map was constructed with the requirement that the graphs $\Gamma(F_n)$
are isomorphic to $\Gamma\big( P_n(f_{-2}) \big)$, where $f_{-2}:z \mapsto
z^2-2$ is the Chebyshev polynomial. The motivation for this example is to
investigate what properties of $Q$-recurrent polynomials will hold for
parameters that imitate the behavior of (non-periodic) postcritically finite
maps. Actually, the construction of $P_{\rm MCheb}$ is very similar to that of
$Q$-recurrent polynomials; the main difference being as follows. Since the
critical orbit of $f_{-2}$ does not return to the center, the first return to
level $n+1$ must be delayed until after the composition of $n$ first return
maps $g_1 \circ \ldots \circ g_n$, when the critical orbit falls in
$Y_0^{(0)}$.

The parameter was chosen so that the first return to level $n+1$ occurs
exactly at this moment; that is, $g_{n+1}$ is precisely $g_1 \circ g_2 \circ
\ldots \circ g_n$ for all $n$. These are the iterates marked with a
check. Moreover, the choice of frame orientations that result in a real
parameter imposes the required label sequence $\big( 2;\, '{{\sf Z}_1}',\,
'{{\sf LZ}_0}',\, '{{\sf LRZ}_0}',\, '{{\sf LRRZ}_0}',\, '{{\sf
LRRRZ}_0}',\,\ldots \big)$. By analogy with the critical orbit of $f_{-2}$,
every nest level of $P_{\rm MCheb}$ has two lateral pieces and the itinerary
of $V_i^n$ includes an infinite number of visits to $V_2^n$ after the first
return to $V_0^{n+1}$.

Since this combinatorial type is admissible, the results of \cite{1st_part}
guarantee an uncountable set of complex parameters with the same
combinatorics. By the rigidity result of Yoccoz, there cannot be other real
polynomials in this class. \\

The methods used to work with $Q$-recurrent polynomials are not enough to
study the nest of $P_{\rm MCheb}$ in a metric sense. In particular, we have
relied on the fact that $K_Q$ had a non-empty interior; this is not the case
for $K_{-2}$. Nevertheless, the analogy is good enough that it is natural to
pose the following.

\begin{conj}
  There are suitable rescalings of the nest pieces of $P_{\rm MCheb}$ such
  that the functions induced from the first return maps converge to $f_{-2}$
  and such that properly rescaled pieces converge to the interval $[-2,2]$ in
  the Hausdorff topology.
\end{conj}

\section{Parameter space}\label{sect:Parameter_Space}
One of the most amazing attributes of complex quadratic dynamics is the
replication of dynamical features in the parameter plane. For instance, the
structure of a limb $L_{p/q}$ reflects the initial steps of the critical orbit
for any parameter contained in it. In \cite{TanLei}, Tan Lei showed that for a
strictly preperiodic parameter $c$, the Julia set of $f_c$ and the Mandelbrot
set exhibit local asymptotic similarity around $c$.

As it has been mentioned, a result of similar nature appears in \cite{W} where
Wenstrom shows that the paranest pieces around the real Fibonacci parameter
$c_{\rm fib}$ are asymptotically similar to the central pieces in the
principal nest of $f_{c_{\rm fib}}$. Thus, $\Delta^n(c_{\rm fib})
\longrightarrow K_{-1}$ in shape and the author exploits this geometric result
to obtain hairiness of $M$ around $c_{\rm fib}$.

This Section discusses a generalization of the above results to the family of
all $Q$-recurrent parameters. Note that the maps $Q$ are dense in $\bdry M$
and that for each one there is an uncountable set of $Q$-recurrent parameters.
 
\begin{thm}\label{thm:Paranest_Shape}
  Let $Q = Q(z)$ be the center of a prime hyperbolic component with critical
  period $m$, and let $c_Q$ be a $Q$-recurrent parameter. Then the paranest
  around $c_Q$ is infinite and the parapieces $\big\{ \Delta^j(c_Q) \big\}$
  converge in shape to the filled Julia set $K_Q$.
\end{thm}

This will require translating the corresponding result obtained in the
dynamical plane to the space of parameters. To do this, we need to introduce
certain auxiliary parapieces; describe in detail the boundary of
$\Delta^j(c_Q)$ and define a map $M_n:\Delta^n(c_Q) \longrightarrow
\mathbb{C}$ that ``rescales'' $\Delta^n$ to a compact set close to $K_Q$.

From the above result follows the possibility of computing the rate of growth
of the paramoduli. Since the paramoduli increase at least linearly, the set of
$Q$-recurrent parameters is a Cantor set of Hausdorff dimension 0.

For the rest of this Section, unless explicitly mentioned, fix a map $Q =
Q(z)$ in the center of a prime hyperbolic component such that the critical
orbit has period $m$; also, $c_Q$ will stand for a fixed $Q$-recurrent
parameter.

\subsection{Auxiliary parapieces}
Consider the first return map $g_{n-1}: V_0^{n-1} \longrightarrow V_0^{n-2}$.
Given that $V_1^n \Subset V_0^{n-1}$, we can study the effect of $g_{n-1}$ on
$V_1^n$; there, the condition of $Q$-recurrency gives:
$$g_{n-1}(V_1^n) \Subset V_2^{n-1}.$$

Applying $g_{n-2}$ we obtain
$$g_{n-2} \circ g_{n-1}(V_1^n) \Subset g_{n-2}(V_2^{n-1}) \Subset V_3^{n-2}.$$

\begin{figure}[h]
  \includegraphics{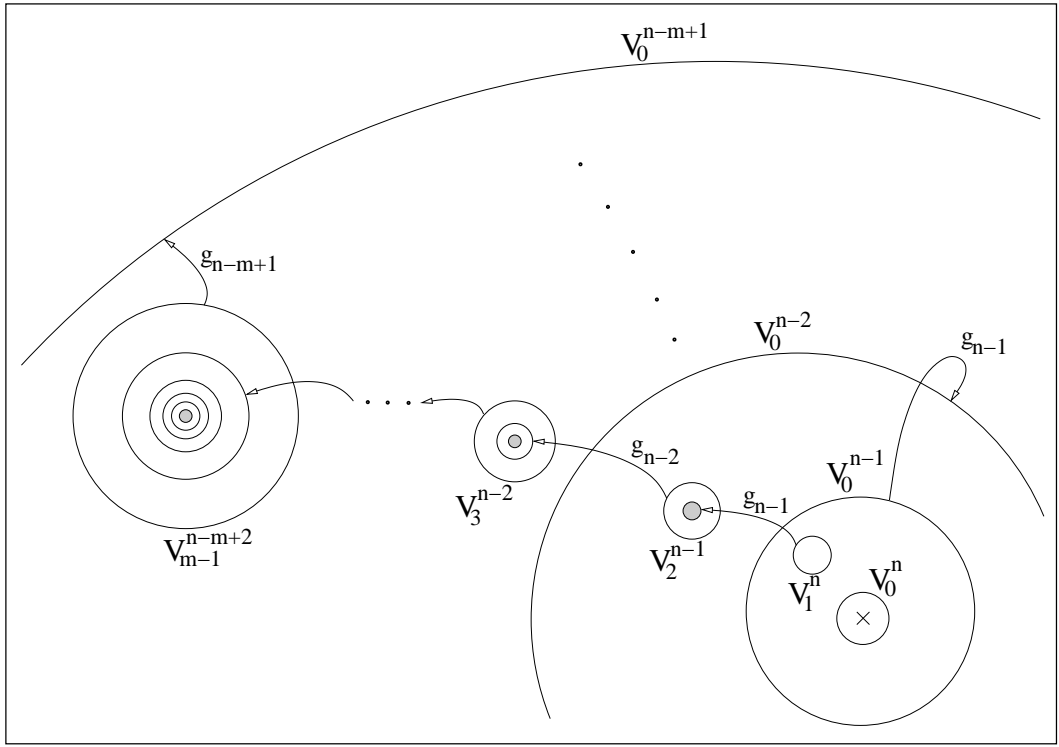}
  \caption[Construction of $U^n$ and $F_{n+2}^*$.]
          {
           \label{fig:Detail_Of_Q-nest} \it Construction of $U^n$ and
           $F_{n+2}^*$. The domain of $g_{n-1}$ is $V^{n-1}$; in particular,
           it maps $V_1^n$ inside $V_2^{n-1}$. Further maps $g_{n-2}, g_{n-3},
           \ldots$ take the current ensemble inside the next piece of previous
           level, until $V_{m-1}^{n-m+2}$. Note that $g_{n-m+1}$ takes
           $V_{m-1}^{n-m+2}$ onto the central piece $V_0^{n-m+1}$, instead of
           a lateral one. This maps the gray piece (the image of $V_1^n$) onto
           $V_0^{n-1}$. Then, we can pull $V_0^n$ back all the way to the
           piece $U^n$ inside $V_1^n$. Also, $U^n$ has a frame $F_{n+2}^*$
           which is the corresponding pull-back of $F_{n+2} \subset
           V_0^n$. Neither $U^n$ nor $F_{n+2}^*$ are drawn. See also Figure
           \ref{fig:Q_Recurrency}.  }
\end{figure}

This procedure can continue further for a total of $m-2$ steps:
$$g_{n-m+2} \circ \ldots g_{n-1}(V_1^n) \Subset V_{m-1}^{n-m+2},$$ where in
fact, the last image is contained in a nest of intermediate pieces inside
$V_{m-1}^{n-m+2}$; see Figure \ref{fig:Detail_Of_Q-nest}. Under $g_{n-m+1}$,
the piece $V_{m-1}^{n-m+2}$ maps onto $V_0^{n-m+1}$ so the intermediate pieces
mentioned will map onto the central pieces $V_0^{n-m+2}, V_0^{n-m+3}, \ldots,
V_0^{n-2}$ and the combined effect on $V_1^n$ will be (recall Definition
\ref{defn:g_{n,k}}):
\begin{equation}\label{g_{n,1}}
  g_{n-m+1} \circ \ldots g_{n-1}(V_1^n) = g_{n,1}(V_1^n) = V_0^{n-1}.
\end{equation}

Since $V_0^n \Subset V_0^{n-1}$, the following is well defined.

\begin{defn}
  Denote by $U^n \Subset V_1^n$ and $F_{n+2}^* \subset U^n$ the
  $(g_{n,1})$-pull-backs of $V_0^n$ and $F_{n+2} \subset V_0^n$, respectively.
  Compare Figure \ref{fig:Detail_Of_Q-nest}.
\end{defn}

Note that $F_{n+2}^*$ is known once the nest structure up to level $n$ is
given. However, assuming that the nest of our parameter displays the
$Q$-recurrency type up to level $n+1$, it is possible to say more. Since
$g_{n+1}(0) = g_{n-m+1} \circ \ldots g_{n-1} \circ g_n (0)$, we must have
\begin{equation}\label{F*}
  g_n(0) \in F_{n+2}^*.
\end{equation}

Let us pass to the parameter plane. Our initial goal is to obtain a precise
control of the combinatorics inside relevant consecutive parapieces. In the
first place, $\Delta^n$ is the set of parameters that have the same nest
combinatorics as $c_Q$, up to the first return $g_n(0)$ to $V_0^{n-1}$.

\begin{defn}
  We introduce two new auxiliary parapieces.
  \begin{itemize}
    \item $\Delta_*^{n+2}$ is the set of parameters such that $g_n(0)$ falls
          inside the frame $F_{n+1} \subset V_0^{n-1}$.
    \item $\Xi^{n+1}$ is the set of parameters such that $g_n(0)$ falls in
          $V_1^n \Subset F_{n+1}$.
  \end{itemize}
\end{defn}

Each region $\Xi^n$ is well defined as a {\bf parapiece} since it represents
the return to an explicit piece of the puzzle. On the other hand,
$\Delta_*^{n+2}$ is actually the union of several parapieces; nevertheless, it
is convenient to regard it as a parapiece to avoid longer descriptions.  With
this in mind, we are interested in the fact that parapieces of consecutive
levels can be described in terms of a single first return map. From
Formulas (\ref{Annuli_Returns}), (\ref{g_{n,1}}) and (\ref{F*}), we obtain:
\begin{equation}\label{eqn:Parameaning}
\begin{array}{lclcl}
  c \in \Delta^n       & \Longleftrightarrow & g_n(0) \in V_0^{n-1} &
                         \Longleftrightarrow & g_{n+1}(0) \in V_0^{n-m} \\
  c \in \Delta_*^{n+2} & \Longleftrightarrow & g_n(0) \in F_{n+1} &
                         \Longleftrightarrow & g_{n+1}(0) \in F_{n-m+2} \\
  c \in \Xi^{n+1}      & \Longleftrightarrow & g_n(0) \in V_1^n &
                         \Longleftrightarrow & g_{n+1}(0) \in V_0^{n-1} \\
  c \in \Delta^{n+1}   & \Longleftrightarrow & g_n(0) \in U^n &
                         \Longleftrightarrow & g_{n+1}(0) \in V_0^n \\
  c \in \Delta_*^{n+3} & \Longleftrightarrow & g_n(0) \in F_{n+2}^* &
                         \Longleftrightarrow & g_{n+1}(0) \in F_{n+2}. \\
\end{array}
\end{equation}

Since
$F_{n+2}^* \subset U^n \Subset V_1^n \Subset F_{n+1} \subset V_0^{n-1}$,
we have the following parapiece inclusions:
\begin{equation}
  \Delta_*^{n+3} \subset \Delta^{n+1} \Subset \Xi^{n+1}
  \Subset \Delta_*^{n+2} \subset \Delta^n.
\end{equation}

\subsection{Shape and paramoduli}
In order to prove Theorem \ref{thm:Paranest_Shape}, let us introduce the map
$M_n: \Xi^{n-1} \longrightarrow \mathbb{C}$, where $\Xi^{n-1}$ belongs to the
paraframe of the fixed parameter $c_Q$. Recall that $\alpha_{n-2} =
\frac{\alpha_Q}{s_{n-2}}$ is the rescaling factor that defines
$\Tilde{F}_n[c_Q] = \alpha_{n-2} \cdot F_n[c_Q]$.

For $c \in \Xi^{n-1}$, let $$M_n(c) = \alpha_{n-2} \cdot g_{n-1}(0)[c].$$

\begin{proof}[Proof of Theorem \ref{thm:Paranest_Shape}.]
  From Table (\ref{eqn:Parameaning}), when $c \in \Delta_*^{n+1}$ the first
  return $g_{n-1}(0)[c]$ is in $F_n$. Recall that for $n$ large,
  $\Tilde{F}_n[c]$ is exponentially close to $K_Q$.

  Fix an $\varepsilon >0$ and find $n$ big enough so that both rescalings
  $\alpha_{n-2}[c] \cdot F_n = \Tilde{F}_n$ and $\alpha_{n-2} \cdot F_n$
  are at most $\frac{\varepsilon}2$-close to each other and to $K_Q$. This
  means that $M_n(\Delta_*^{n+1})$ is a compact set $\varepsilon$-close to
  $K_Q$. \\

  By definition, the parapiece $\Xi^{n-1}$ is the set of parameters $c$ for
  which $g_{n-2}(0)[c]$ falls on the lateral piece $V_1^{n-2}$. Since this map
  is a first return, Proposition \ref{prop:Identify_Bdries} implies that the
  correspondence $c \mapsto g_{n-2}(0)[c]$ is univalent in
  $\Xi^{n-1}$. Moreover, $V_1^{n-2}[c]$ is at a definite distance away from
  the central piece $V_0^{n-2}[c]$ for all $c$, so the image of $\Xi^{n-1}$
  under $c \mapsto g_{n-2}(0)[c]$ is uniformly far from 0 and similarly for
  all further iterations up to the first return $g_{n-1}(0)[c]$. Again, we can
  use Proposition \ref{prop:Identify_Bdries} to deduce that $M_n$ is
  univalent in its entire domain.

  Since $n$ is big, the modulus of the annulus $\big( {\rm int} V_0^{n-3}
  \setminus F_n \big)[c]$ is large for every $c \in \Xi^{n-1}$. In
  particular, since $\Tilde{F}_n[c]$ has bounded diameter, this implies that
  the distance between the point $\alpha_{n-2}[c] \cdot g_{n-1}(0)[c] \in
  \bdry \Tilde{V}_0^{n-3}[c]$ and the curve $\bdry \Tilde{F}_n[c]$ is
  exponentially big for $c \in \bdry \Xi^{n-1}$. Let $d_n$ be the minimum of
  these distances over all $c$. Then, $M_n(\bdry \Xi^{n-1})$ and $M_n(\bdry
  \Delta_*^{n+1})$ are at least a distance $(d_n - \varepsilon) \sim d_n
  \nearrow \infty$ apart.  We can conclude that the modulus of $M_n\big( {\rm
  int\,} \Xi^{n-1} \setminus \Delta_*^{n+1} \big)$ is arbitrarily large and so
  will be the modulus of ${\rm int\,} \Xi^{n-1} \setminus \Delta_*^{n+1}$. \\

  We have shown that the map $M_n$ is univalent in its domain and the modulus
  of ${\rm int\,} \Xi^{n-1} \setminus \Delta_*^{n+1}$ is big. Then, by the
  Koebe distortion Theorem, $M_n$ is asymptotically linear in a neighborhood
  of $ \Delta_*^{n+1}$. Since $M_n\big( \Delta_*^{n+1} \big)$ is
  $\varepsilon$-close to $K_Q$, we can conclude the proof.
\end{proof}

As an immediate consequence of this control over the shape of parapieces, we
can compute the rate of growth of principal paramoduli $\mu_n$.

\begin{corol}\label{corol:Para_Growth}
  The annuli of consecutive parapieces in the nest of a $(z^2-1)$-recurrent
  map grow linearly at the rate
  $$\lim_{n \rightarrow \infty} \frac{\mu_n}n = 2\frac{\ln 2}3.$$

  For any other $Q$-recurrent map (where $Q$ has critical orbit of period $m
  \geq 3$) the moduli grow exponentially at the rate
  $$\lim_{n \rightarrow \infty} \frac{\mu_n}{\mu_{n+1}} = \kappa_m,$$
  where $\kappa_m$ is the same constant as in Theorem \ref{thm:Moduli_Growth},
  converging to $\frac32$ as the period $m$ of $Q$ increases.
\end{corol}

\begin{proof}
  First note that, although $U^n$ is defined as a pull-back of $V_0^n$,
  relation \ref{F*} shows that this piece is just $g_n(V_0^{n+1})$.

  Now, when $c \in \Delta^n$, the first return $g_n(0)$ falls in $V_0^{n-1}$.
  For $c \in \Delta^{n+1}$, $g_n(0)$ is in $U^n$. From the previous result, $c
  \mapsto g_n(0)[c]$ is an almost linear map taking the annulus $\big(
  \Delta^n \setminus \Delta^{n+1} \big)$ close to $\big( V_0^{n-1} \setminus
  U^n \big)[c_Q]$. Therefore
  $$
    {\rm mod} \big( \Delta^n \setminus \Delta^{n+1} \big) \sim
    {\rm mod} \big( V_0^{n-1} \setminus U^n \big) \sim
    2\,{\rm mod} \big( V_0^n \setminus V_0^{n+1} \big).
  $$

  Corollary \ref{corol:Para_Growth} now follows from Theorem
  \ref{thm:Moduli_Growth}.
\end{proof}

\subsection{Auto-similarity in the Mandelbrot set}\label{sect:Auto_Similarity}
The discovery that parapieces around $c_{\rm fib}$ are similar to the Julia
set of $-1$ revealed one more level of complexity in the structure of $M$
since it relates the dynamics of two different parameters. In this Subsection
we use our results to take one further step. Having at our disposal an
infinite collection of superattracting parameters, we reveal an interesting
relation between two arbitrary parameters on $\bdry M$ whose combinatorics can
be completely dissimilar.

The $Q$-recurrency phenomenon is not restricted to the Cantor sets described
so far. As part of the proof of the next Theorem, we will show that parapieces
whose shape approximates $K_Q$ are dense on $\bdry M$. This requires relaxing
the definition of $Q$-recurrency which assumes that the correct combinatorics
start from level 0. Instead, we allow critical orbits that behave arbitrarily
for several levels before settling in the desired $Q$-recurrent pattern. This
critical behavior is referred to as {\it generalized $Q$-recurrency}. The
assertion of Theorem \ref{thm:M_Similarity} follows; it can be interpreted as
saying that the geometry of most Julia sets is replicated near arbitrary
locations of the boundary of $M$.

\begin{thm}\label{thm:M_Similarity}
  Let $c_1, c_2 \in \bdry M$ be two parameters such that $f_{c_2}$ has no
  indifferent periodic orbits that are rational or linearizable. Then there
  exists a sequence of parapieces $\{ \Upsilon_1, \Upsilon_2, \ldots \}$ (most
  likely not nested) converging to $c_1$ as compact sets, but such that
  $\Upsilon_n \longrightarrow K_{c_2}$ in shape.
\end{thm}

\begin{proof}
  It is not difficult to obtain the result of Theorem
  \ref{thm:Convergence_Of_G_n} in more generality. In fact, inside any ball
  $B_{\varepsilon}(c)$ with $c \in \bdry M$, we can find a system of
  parameters for which the first return maps converge (after scaling) to a
  given superattracting map $Q$.

  To see this, simply consider a tuned copy of $M$ contained in
  $B_{\varepsilon}$. All parameters in this copy $M'$, are renormalizable by
  the same combinatorics. In particular, there will be parameters whose
  renormalization is hybrid equivalent to a $Q$-recurrent map. For these
  parameters a high level of the frame will contain a substructure whose graph
  is isomorphic to $\Gamma_0(Q)$ and we can start the same construction as in
  the proof of Theorem \ref{thm:Convergence_Of_G_n} to produce frame-like
  structures whose graphs are isomorphic to $\Gamma_n(Q)$. Since the
  combinatorics is prescribed by a polynomial, there can be no obstructions
  just as in the original case. Then, the rescaled first return maps
  will converge to $Q$ as before. Moreover, we can translate the shape
  property to the parameter plane.

  Note that this argument is equivalent to prescribing the itineraries of nest
  pieces arbitrarily on the initial levels and then proving that they can be
  admissibly extended on subsequent levels to match the pattern given in
  Formula (\ref{eqn:Q_Itineraries}). \\

  Now consider the filled Julia set of $f_{c_2}$. We know from
  \cite{Julia_set_continuity} that there is a sequence $\{ Q_1, Q_2, \ldots
  \}$ of superattracting polynomials in a prime hyperbolic component of $M$
  such that $K_{c_2}$ can be arbitrarily approximated by filled Julia sets:
  $K_{Q_n} \longrightarrow K_{c_2}$. To fix ideas, let us choose subindices so
  that the Hausdorff distance is ${\rm dist}_H(K_{Q_n},K_{c_2}) < \frac1{2n}$.

  For any $n$, consider the ball $B_{\frac1n}(c_1)$ and locate a generalized
  $Q_n$-recurrent parameter $s_n$. By going to a deep enough level, we can
  find some parapiece $\Upsilon_n \subset B_{\frac1n}$ around $s_n$ whose
  shape is $\big( \frac1{2n} \big)$-close to $K_{Q_n}$; that is, so that there
  is a rescaling $\Tilde{\Upsilon}_n$ of $\Upsilon_n$ for which ${\rm
  dist}_H\big( \Tilde{\Upsilon}_n, K_{Q_n} \big) < \frac1{2n}$.

  Since $\Upsilon_n \subset B_{\frac1n}$, the sequence $\{ \Upsilon_n \}$
  consists of parapieces that get arbitrarily small and converge to $c_1$,
  while at the same time ${\rm dist}_H\big( \Tilde{\Upsilon}_n,K_{c_2} \big) <
  \frac1n$, so $\Upsilon_n \longrightarrow K_{c_2}$ in shape.
\end{proof}

\section*{Appendix}
\setcounter{subsection}{0} \setcounter{thm}{0}
\renewcommand{\thesubsection}{\Alph{subsection}}
\renewcommand{\thethm}{\thesubsection.\arabic{thm}}
The theorems of previous sections rely on several notions and results of
complex analysis. In this Appendix we will describe the necessary ideas on
which the text relies. For references and proofs of these results, the reader
can consult \cite{A_qc_book}, \cite{DH_thurston} and \cite{LV_book}.

\subsection{Carath\'eodory topology}
\setcounter{thm}{0} A sequence of pointed disks $\big\{ (U_n,x_n) \big\}$ is
said to converge to the pointed disk $(U,x)$ in the Carath\'eodory topology
if: \renewcommand{\theenumi}{\alph{enumi}}
\begin{enumerate}
  \item $x_n \longrightarrow x$.
  \item For every compact $K \subset U$ there is an $N$ such that $K \subset
        U_n$ for $n>N$.
  \item If $ V \ni x$ is an open connected region and $V \subset U_n$ for
	infinitely many $n$, then $V \subset U$.
\end{enumerate}

The interpretation of convergence in this topology is as follows. Consider the
complements $\hat{\mathbb{C}} \setminus U_n$ which are compact sets converging
to $X$ in the Hausdorff topology. To satisfy the above conditions, $x \nin X$
and $U$ is the component of $\hat{\mathbb{C}} \setminus X$ that contains
$x$. \\

We can describe a similar convergence for a sequence of annuli $\{A_n\}$. In
this case, we require that $\hat{\mathbb{C}} \setminus A_n \longrightarrow Y$
in the Hausdorff sense and that the core geodesics $\gamma_n$ of the annuli
$A_n$ converge to a non-degenerate loop $\gamma \subset Y$. The Carath\'eodory
limit will be the doubly connected component of $\hat{\mathbb{C}} \setminus Y$
containing $\gamma$.

\subsection{Modulus and capacity as conformal invariants}
\setcounter{thm}{0}

Let $R \subset \mathbb{C}$ be a doubly connected region in the plane. Consider
the family $\Gamma$ of curves $\gamma \subset R$ whose endpoints are on the
boundary of $R$. Given a conformal metric $\rho$ on $R$, we can define the
length of $\Gamma$ as $L_{\rho}(\Gamma) = \inf_{\gamma \in \Gamma}
\int_{\gamma}\rho|dz|$ and the area of $R$ as $A_{\rho}(R) = \int \int_{A}
\rho^2 dxdy$.

\begin{defn}
  The modulus of $R$ is
  $$
    {\rm mod}(R) =
    \sup_{\rho} \frac{\big( L_{\rho}(\Gamma) \big)^2}{A_{\rho}(R)}
  $$
  where the supremum is taken over all conformal metrics $\rho$ with
  non-degenerate area: $A_{\rho}(R) \neq 0,\infty$.
\end{defn}

It can be shown that the modulus is a conformal invariant. As a consequence,
we can give an alternative definition as follows. Consider the Riemann map
$\varphi:R \longrightarrow S$ where $S$ is the round annulus $S = \{ 1< z < r
\}$. Then ${\rm mod(R)} = 2\pi \ln r$. \\

Another conformal invariant, this time of topological disks, is the {\it
capacity} or {\it conformal radius}.

\begin{defn}
  Let $U \in \mathbb{C}$ be a simply connected domain, $z \in U$ and
  $\varphi:(\mathbb{D},0) \longrightarrow (U,z)$ the Riemann map from the unit
  disk to $U$ that satisfies $\varphi'(0) > 0$. The {\bf capacity} of $U$ with
  respect to $z$ is
  $${\rm cap}_z(U) = \ln \varphi'(0).$$
\end{defn}

Of interest to us, will be the following property of capacities and moduli.

\begin{thm}
  Both capacity and modulus are quantities that vary continuously with respect
  to the Carat\'eodory topology.
\end{thm}

\subsection{Gr\"otzsch inequality}
\setcounter{thm}{0}

The following result (and its quantitative version) is essential to estimate
the modulus of an annulus that is split in subannuli.

\begin{thm}{\bf Extended Gr\"otzsch Inequality:}
  {
  \renewcommand{\theenumi}{\alph{enumi}} Let $K \subset \mathbb{C}$ be a
   simply connected compact set and denote ${\rm int}_0 K$ the component
   of its interior that contains 0.
  \begin{enumerate} 
    \item \label{mod_ineq} Consider two topological disks
          $0 \in U \subset {\rm int}_0 K \subset V$. Then
          $$
            {\rm mod} (V \setminus U) \geq
            {\rm mod} (K \setminus U) +
            {\rm mod} (V \setminus K).
          $$
    \item Let $\{U_n\}$ and $\{V_n\}$ be two sequences of nested topological
          disks satisfying
          \begin{itemize}
            \item $0 \in U_n \subset {\rm int}_0 K$ and
                  ${\rm diam}\, U_n \searrow 0$
            \item $K \subset V_n$ and
                  ${\rm dist}(K, \bdry V_n) \nearrow \infty$.
          \end{itemize}
          Then the deficit in the Gr\"otzsch inequality tends to
          \begin{eqnarray*}
            \lim_{n \rightarrow \infty} \Big(
            {\rm mod} (V_n \setminus U_n) -
            {\rm mod} (K \setminus U_n) -
            {\rm mod} (V_n \setminus K) \Big) = \\ \\
            |{\rm cap}_0 ({\rm int}_0 K)| +
            |{\rm cap}_{\infty} (K)|.
          \end{eqnarray*}
  \end{enumerate}
  }
\end{thm}

An important observation is the fact that equality in (\ref{mod_ineq}) is
achieved if and only if $\bdry K \subset V \setminus U$ maps to a centered
circle under the Riemann map.

\subsection{Koebe distortion Theorem}
\setcounter{thm}{0}

\begin{defn}
  Given an analytic univalent map $\varphi$ between regions $U$ and $V$, the
  {\bf distortion} of $\varphi_0$ is defined as:
  $$
    {\rm Dist}(\varphi_0) =
    \sup_{x,y \in U} \frac{\varphi_0'(x)}{\varphi_0'(y)}.
  $$
\end{defn}

Koebe's Theorem provides great control of the distortion when there is enough
space between $U$ and $V$.

\begin{thm}
  Let $U$ and $V$ be two topological disks with $U \setminus V$. Then there is
  a constant $C$ such that for any univalent map $\varphi(U) = V$,
  $$
    {\rm Dist}(\varphi) < C.
  $$ Moreover, $C = 1 + O\big( e^{-{\rm mod}(V \setminus U)} \big)$ as the
  modulus goes to $\infty$.
\end{thm}

\subsection{Teichm\"uller space}
\setcounter{thm}{0}

The Teichm\"uller space of a Riemann surface carries a great deal of
structural information. Here we focus in the case that the surface $S$ is the
complex plane punctured at a finite set $\mathcal{O}$. Then, the Teichm\"uller
space $\mathcal{T}_S$ can be described as a quotient of the space of
quasiconformal deformations of $S$ (i.e. the family of maps $\{h:S
\longrightarrow \mathbb{C}\,|\, h \text{ is a qc homeomorphism} \}$), where
two deformations $h_1$ and $h_2$ are identified if and only if there is a
conformal change of coordinates $\varphi:\mathbb{C} \longrightarrow
\mathbb{C}$ such that $\varphi \circ h_1$ is isotopic to $h_2$ relative to the
puncture set $h_2(\mathcal{O})$.

\note{} The coordinate changes $\varphi$ are affine maps, so the deformation
  of $\mathcal{O}$ within a class is determined up to translation and complex
  scaling. Therefore, we can normalize a deformation $h$ by requiring that $h$
  fixes two distinguished points in $\mathcal{O}$. These could be, for
  instance, the critical point and critical value of $Q$ in the case that
  $\mathcal{O}$ is the postcritical set of a hyperbolic map $Q$. \\

It is fundamental to consider an alternate description of $\mathcal{T}_S$ in
terms of Beltrami differentials. Fix two almost complex structures on $S$
determined by their Beltrami coefficients $\mu \frac{d\overline{z}}{dz}$ and
$\nu \frac{d\overline{z}}{dz}$. Assume that they are related by $\nu =
{\overline{h}}^{\,*}\mu$, where $\overline{h}:S \longrightarrow S$ is a
quasiconformal self homeomorphism of $S$ which is homotopic to id relative to
$\mathcal{O}$. Then, the straightening maps $h_{\mu}$ and $h_{\nu}$ are two
quasiconformal deformations of $S$ in the same equivalence class in
$\mathcal{T}_S$.

Conversely, we can associate to a deformation $h$ the almost complex structure
$h^*\sigma = \frac{\overline{\partial}h}{\partial h}\frac{d\overline{z}}{dz}$
where $\sigma$ is the standard structure. It is easy to verify that this
correspondence lifts to the equivalence classes where it induces a bijection.

\end{document}